\documentclass{article}

\usepackage{amssymb}
\usepackage{amsmath}
\usepackage[all]{xy}
\usepackage[frenchb]{babel}
\usepackage{a4wide}
\usepackage{enumerate,fancyhdr}

\begin{document}

\title{Torsion pour les vari\'et\'es ab\'eliennes de type I et II}

\author{Marc Hindry\footnote{marc.hindry@imj-prg.fr}, Nicolas Ratazzi \footnote{nicolas.ratazzi@math.u-psud.fr}
}

\renewcommand{\epsilon}{\varepsilon}
\renewcommand{\tilde}{\widetilde}

\newcommand{\Lie}{\textnormal{Lie}}
\newcommand{\Gal}{\textnormal{Gal}}
\newcommand{\tors}{\textnormal{tors}}
\newcommand{\Ker}{\textnormal{Ker}}
\newcommand{\red}{\textnormal{r\'{e}d}}
\newcommand{\Z}{\mathbb{Z}}
\newcommand{\Q}{\mathbb{Q}}
\newcommand{\C}{\mathbb{C}}
\newcommand{\R}{\mathbb{R}}
\newcommand{\N}{\mathbb{N}}
\renewcommand{\P}{\mathbb{P}}
\newcommand{\F}{\mathbb{F}}
\newcommand{\G}{\mathbb{G}}
\newcommand{\SD}{\mathbb{S}}
\renewcommand{\hat}{\widehat}

\newcommand{\T}{\mathcal{T}}

\newcommand{\K}{\overline{K}}

\newcommand{\mP}{\mathcal{P}}
\newcommand{\mA}{\mathcal{A}}
\renewcommand{\O}{\mathcal{O}}
\renewcommand{\L}{\mathcal{L}}
\newcommand{\V}{\mathcal{V}}
\newcommand{\kk}{\mathbb{K}}

\newcommand{\im}{\textnormal{Im}}

\newcommand{\dimrel}{\textnormal{dimrel}}
\newcommand{\m}{\mathfrak{m}}
\newcommand{\I}{\mathfrak{I}}

\newcommand{\hl}{\mathfrak{h}_{\ell}}
\newcommand{\el}{\mathfrak{e}_{\ell}}
\newcommand{\gsp}{\mathfrak{gsp}}
\newcommand{\n}{\mathfrak{n}}
\newcommand{\Tr}{\textnormal{Tr}}
\newcommand{\ssp}{\mathfrak{sp}}
\newcommand{\ssl}{\mathfrak{sl}}

\newcommand{\VV}{\mathsf{V}}

\newcommand{\mmid}{\|\|}

\newcommand{\Proj}{\textnormal{Proj\,}}
\newcommand{\spec}{\textnormal{spec\,}}
\newcommand{\Deg}{\textnormal{Deg\,}}
\newcommand{\e}{\varepsilon}
\newcommand{\ab}{\textnormal{ab}}
\newcommand{\End}{\textnormal{End}}
\renewcommand{\H}{\textnormal{H}}
\newcommand{\MT}{\textnormal{MT}}
\newcommand{\codim}{\textnormal{codim\,}}
\newcommand{\Hom}{\textnormal{Hom}}
\newcommand{\Aut}{\textnormal{Aut}}
\newcommand{\Hdg}{\textnormal{Hdg}}
\newcommand{\GL}{\textnormal{GL}}
\newcommand{\SL}{\textnormal{SL}}
\newcommand{\SLG}{\textnormal{SLG}}
\newcommand{\Sp}{\textnormal{Sp}}
\newcommand{\GSp}{\textnormal{GSp}}
\newcommand{\Res}{\textnormal{Res}}
\newcommand{\tr}{\textnormal{Tr}}
\newcommand{\frob}{\textnormal{Frob}}
\renewcommand{\ss}{\sigma}
\newcommand{\pr}{\textnormal{pr}}
\newcommand{\hdg}{\textnormal{Hdg}}
\newcommand{\mult}{\textnormal{mult}}

\newcommand{\rang}{\textnormal{rang\,}}

\newcommand{\W}{\mathcal{W}}

\renewcommand{\text}{\textnormal}

\newcommand{\ee}{\overline{e}}
\newcommand{\ff}{\overline{f}}
\newcommand{\tf}{\hat{f}}
\newcommand{\te}{\hat{e}}

\newcommand{\cl}{\mathcal{L}}

\newcommand{\II}{\textnormal{II}}
\newcommand{\typeI}{\textnormal{I}}

\newcommand{\que}{\asymp}
\newcommand{\pre}{\circeq}
\newcommand{\qui}{\lesssim}

\newtheorem{theo}{Th{\'e}or{\`e}me} [section]
\newtheorem{lemme}[theo]{Lemme}
\newtheorem{conj}{Conjecture}[section]
\newtheorem{prop}[theo]{Proposition}
\newtheorem{cor}[theo]{Corollaire}
\newtheorem{cons}[theo]{Cons\'equence}

\newcommand{\defi}{\addtocounter{theo}{1}{\noindent \textbf{D{\'e}finition \thetheo\ }}}
\newcommand{\rem}{\addtocounter{theo}{1}{\noindent \textbf{Remarque \thetheo\ }}}
\newcommand{\exemple}{\addtocounter{theo}{1}{\noindent \textbf{Exemple \thetheo\ }}}

\newcommand{\demo}{\noindent \textit{D{\'e}monstration} : }
\newcommand{\findemo}{\hfill \Box}

\maketitle

\hrulefill

\noindent \textbf{R\'esum\'e : } 
Soit $A$ une vari\'et\'e ab\'elienne d\'efinie sur un corps de nombres $K$, le nombre de points de torsion
d\'efinis sur une extension finie $L$ est born\'e polynomialement en terme du degr\'e $[L:K]$ de $L$ sur $K$. Sous les trois hypoth\`eses suivantes nous calculons l'exposant optimal dans cette borne, en terme de la dimension des sous-vari\'et\'es ab\'eliennes de $A$ et de leurs anneaux d'endomorphismes. Les trois hypoth\`eses faites sur A sont les suivantes : (1) $A$ est g\'eom\'etriquement isog\`ene \`a un produit de vari\'et\'es ab\'eliennes simples de type I ou II dans la classification d'Albert ; (2) $A$ est de ``type Lefschetz" c'est-\`a-dire que le groupe de Mumford-Tate est le groupe des similitudes symplectiques commutant aux endomorphismes ; (3) $A$ v\'erifie la conjecture de Mumford-Tate. Le r\'esultat est notamment inconditionnel (i.e. ces trois hypoth\`eses sont v\'erifi\'ees) pour un produit de vari\'et\'es ab\'eliennes simples de type I ou II et de dimension relative impaire. Par ailleurs nous prouvons, en \'etendant des r\'esultats de  Serre, Pink et Hall,  la  conjecture de Mumford-Tate pour quelques nouveaux cas de vari\'et\'es ab\'eliennes de type Lefschetz.
  
\bigskip

\hrulefill

\bigskip

\noindent \textbf{Abstract : }
Let $A$ be an abelian variety defined over a number field $K$, the number of torsion points rational over a finite extension $L$ is bounded polynomially in terms of 
the degree $[L:K]$ of $L$ over $K$.  Under the following three conditions we compute the optimal exponent for this bound, in terms of the dimension of abelian subvarieties and their endomorphism rings. The three hypotheses are the following : (1) $A$ is geometrically isogenous to a product of simple abelian varieties of type I or II, according to Albert classification; (2) $A$ is of ``Lefschetz type", that is, the Mumford-Tate group is the group of symplectic similitudes which commute with the endomorphism ring. (3) $A$ satisfies the Mumford-Tate conjecture. This result is unconditional (i.e. the three hypotheses are satisfied) for a product of simple abelian varieties of type I or II with odd relative dimension. Further, building on work of Serre, Pink and Hall,  we also prove    the Mumford-Tate conjecture  for a few new cases of abelian varieties of Lefschetz type.

\hrulefill

\section{Introduction}

\noindent Soit $A/K$ une vari\'et\'e ab\'elienne, d\'efinie sur un corps de nombres $K$, de dimension $g\geq 1$. Le classique th\'eor\`eme de Mordell-Weil assure que le groupe $A(F)$ des points $F$-rationnels de $A$ est de type fini pour toute extension finie $F/K$. Un probl\`eme naturel qui se pose est de d\'eterminer le sous-groupe de torsion $A(F)_{\tors}$. Un premier probl\`eme consiste en fait \`a borner explicitement le cardinal de $A(F)_{\tors}$ lorsque $A$ ou $F$ varient. Comme dans les articles \cite{ratazziens}, \cite{HR} et \cite{hindry-ratazzi1} auxquels ce papier fait suite, nous nous int\'eressons ici au cas o\`u l'on fixe $A$ et o\`u l'on fait varier $F$ parmi les extensions finies de $K$ ; l'objectif \'etant d'obtenir une borne  avec une d\'ependance explicite et optimale  en le degr\'e $[F:K]$. Introduisons maintenant l'invariant que nous allons \'etudier.

\medskip

\defi On pose 
\[\gamma(A)=\inf\left\lbrace x>0\, | \,  \forall F/K \text{ finie, }\ \left|A(F)_{\tors}\right|\ll [F:K]^x\right\rbrace.\]

\noindent  La notation $\ll$ signifie qu'il existe  une constante  $C$, ne d\'ependant que de $A/K$, telle que l'on a $\left|A(F)_{\tors}\right|\leq C [F:K]^x$. 

\medskip

\noindent On peut traduire la d\'efinition en le fait que $\gamma(A)$ est l'exposant le plus petit possible tel que pour tout $\epsilon>0$, il existe une constante $C(\epsilon)=C(A/K,\epsilon)$ telle que pour toute extension finie $F/K$ on a
\[\left|A(F)_{\tors}\right|\leq C(\epsilon) [F:K]^{\gamma(A)+\epsilon}.\]

\noindent Un r\'esultat g\'en\'eral d\^u \`a Masser \cite{lettremasser} donne une borne simple~:
\[\gamma(A)\leq\dim A.\]
\noindent Cette borne est optimale lorsque $A$ est une puissance d'une courbe elliptique avec multiplication complexe; il est fort probable que la borne de Masser n'est jamais optimale dans les autres cas. L'invariant $\gamma(A)$ est calcul\'e dans \cite{hindry-ratazzi1} pour un produit de courbes elliptiques et, de mani\`ere diff\'erente, dans \cite{ratazziens} pour une vari\'et\'e ab\'elienne de type CM et dans \cite{HR} pour une vari\'et\'e ab\'elienne ``g\' en\'erique". Le probl\`eme analogue pour les  modules de Drinfeld est trait\'e dans \cite{florian}. 

\medskip

\noindent Nous notons $\End(A)$ l'anneau d'endomorphismes de $A/\bar K$ et $\End^0(A):=\End(A)\otimes\Q$.  On \'etudie dans ce texte les deux  cas suivants:
\begin{itemize}
\item La vari\'et\'e ab\'elienne  $A$ est de type I,  ie  $E:=\End(A)\otimes \Q$ est un corps totalement r\'eel de degr\'e $e$; en particulier $g=he$ o\`u $h$ est entier. 
\item La vari\'et\'e ab\'elienne  $A$ est de type II,  ie   $D:=\End(A)\otimes \Q$ est  une alg\`ebre de quaternions totalement ind\'etermin\'ee de centre un  corps totalement r\'eel $E$ de degr\'e $e$; en particulier $g=2he$ o\`u $h$ est entier. 
\end{itemize}

\medskip

\noindent Les deux autres cas (type III, ie alg\`ebre de quaternions totalement d\'efinie, et type IV, ie alg\`ebre \`a division sur un corps CM) sont de nature diff\'erente et seront trait\'es dans une autre publication.

\medskip

\defi L'entier $h$ pr\'ec\'edent s'appelle la {\it dimension relative} de $A$. 

\medskip

\noindent On sait que, pour une telle vari\'et\'e ab\'elienne de type I ou II, le groupe de Hodge est toujours contenu dans $\Res_{E/\Q}{\Sp_{2h}}$ et est g\'en\'eriquement \'egal \`a ce dernier groupe. Lorsque $h$ est {\it impair\/}  ou \'egal \`a 2, on sait  que le groupe de Hodge est effectivement  $\Res_{E/\Q}{\Sp_{2h}}$ et, de plus, que la conjecture de Mumford-Tate est vraie (cf. le th\'eor\`eme \ref{mt12} plus loin). Notons que le cas de type I avec $e=1$ correspond aux vari\'et\'es ab\'eliennes de type GSp (dans la terminologie de \cite{HR}).

\medskip

\noindent Pour tout premier $\ell$, on note  $T_{\ell}(A)=\varprojlim A[\ell^n]$ le module de Tate $\ell$-adique de $A$; consid\'erant l'action du groupe de Galois $G_K :=\Gal(\bar K/K)$ sur les points de $\ell^{\infty}$-torsion, on associe naturellement \`a $A/K$, la repr\'esentation $\ell$-adique 
\[\rho_{\ell^{\infty} ,A}: G_K \rightarrow\GL(T_{\ell}(A))\simeq\GL_{2g}(\Z_{\ell})\]
\noindent On pose  $V_{\ell}(A)=T_{\ell}(A)\otimes_{\Z_{\ell}}\Q_{\ell}$. Rappelons que $V_{\ell}(A)$ est naturellement le dual de $H^1_{\rm\acute{e}t}(A\times\bar{K},\Q_{\ell})$ et que, si l'on note
$V=V(A)=H_1(A,\Q)$, alors, comme $\Q$-espace vectoriel, $V\cong \Q^{2g}$, et $V\otimes_{\Q}\Q_{\ell}$ s'identifie naturellement \`a $V_{\ell}(A)$. Si $\psi:V\times V\rightarrow\Q$ est la forme symplectique associ\'ee \`a une polarisation, sa tensoris\'ee par $\Q_{\ell}$ s'identifie avec la forme symplectique provenant de l'accouplement de Weil associ\'e \`a la m\^eme polarisation $T_{\ell}(A)\times T_{\ell}(A)\rightarrow \Z_{\ell}(1)\cong\Z_{\ell}$.

\medskip

\defi Le {\it groupe de Lefschetz} de $A$, not\'e $L(A)$, est le commutateur\footnote{Ce groupe alg\'ebrique est connexe, sauf pour les vari\'et\'es de type III, Cf Murty \cite{murty1} et Milne \cite{milne99}.}, en tant que $\Q$-groupe alg\'ebrique, de $\End^0(A)$ dans $\Sp(V,\psi)=\Sp_{2g}$.

\medskip

\noindent On note \'egalement $G_{\ell}$ ou $G_{\ell,A}$ si besoin, l'enveloppe alg\'ebrique de l'image de $\rho_{\ell^{\infty}}$,  et $H_{\ell,A}$   la composante neutre de $G_{\ell}\cap\SL(V_{\ell})$, c'est-\`a-dire :
\[G_{\ell,A}:=\left(\overline{\rho_{\ell,A}\left(G_K\right)}^{\text{Zar}}\right)^0\qquad{\rm et}\qquad H_{\ell,A}=\left(G_{\ell}\cap\SL(V_{\ell})\right)^0.\]

\noindent On note $\MT(A)$ le groupe de Mumford-Tate et $\Hdg(A)$ le groupe de Hodge d'une vari\'et\'e ab\'elienne $A$ (cf. par exemple \cite{pink} et \cite{moonen}). La relation entre les deux groupes est $\MT(A)=\G_m\cdot\Hdg(A)$ et $\Hdg(A)$ est la composante neutre de $\MT(A)\cap\SL(V)$.

\medskip

\noindent Le groupe des homoth\'eties est contenu dans $G_{\ell,A}$; une autre mani\`ere de d\'efinir le $\Q_{\ell}$-groupe alg\'ebrique $H_{\ell,A}$, pour se ramener dans $\Hdg(A)$, est de le voir comme l'enveloppe alg\'ebrique de $\rho(\Gal(\bar{K}/K(\mu_{\ell^{\infty}}))$, c'est-\`a-dire la composante neutre de la cl\^oture de Zariski de l'image de Galois dans $\GL_{2g,\Q_{\ell}}$. On sait alors que l'on a toujours
\begin{equation*}
H_{\ell,A}\subset\Hdg(A)_{\Q_{\ell}}\subset L(A)_{\Q_{\ell}}.
\end{equation*}

\noindent Remarquons que ces deux inclusions n'ont pas le m\^eme statut : la conjecture  de Mumford-Tate pr\'edit que la premi\`ere inclusion est toujours une \'egalit\'e, alors que la seconde inclusion  peut ne pas \^etre une \'egalit\'e. 

\medskip

\defi Nous dirons que $A$ est de {\it  type   Lefschetz} si $\Hdg(A)= L(A)$ et {\it pleinement de type Lefschetz} si pour tout premier $\ell$, on a $H_{\ell,A}=L(A)_{\Q_{\ell}}$.

\medskip

\rem Dans la d\'efinition pr\'ec\'edente, il suffit en fait de supposer que $H_{\ell,A}=L(A)_{\Q_{\ell}}$ pour un premier $\ell$ (cf. theorem 4.3 de \cite{lp}). 

\medskip

\noindent Ainsi une vari\'et\'e ab\'elienne de type Lefschetz est pleinement de type Lefschetz si elle v\'erifie la conjecture de Mumford-Tate.

\begin{theo}\label{maintheo} Soit  $A$ une vari\'et\'e ab\'elienne g\'eom\'etriquement simple de type I ou II d\'efinie sur un corps de nombres, dont le centre de l'alg\`ebre d'endomorphismes est un corps de nombres totalement r\'eel $E$ de degr\'e $e=[E:\Q]$, posons $d:=\sqrt{[\End^0(A):E]}$ (ie $d=1$ ou 2). On suppose que $A$ est pleinement de type Lefschetz, on a alors 
\begin{equation*}
\gamma(A)=\frac{2dhe}{1+2eh^2+he}=\frac{2\dim A}{\dim\MT(A)}.
\end{equation*}
\end{theo}

\medskip

\rem Il est int\'eressant d'observer que l'exposant $\gamma(A)$ est beaucoup plus petit que la borne $\gamma(A)\leq g$ donn\'ee par Masser ; par exemple, pour toute vari\'et\'e ab\'elienne v\'erifiant les hypoth\`eses du th\'eor\`eme, on a pour le type I (resp. le type II) $\gamma(A)<2/3$ (resp. $\gamma(A)<4/3$).

\medskip

\noindent Pour \'enoncer la deuxi\`eme partie du corollaire suivant, nous introduisons la notation~: 
\begin{equation}\label{defdeS}
\Sigma=\left\{g\geq 1\;|\; \exists k\geq 3,\;{\rm impair},\; \exists a\geq 1,\; 2g=(2a)^{k} \text{ ou } 2g={2k\choose k}\right\}.
\end{equation}
\noindent Pour \'enoncer la troisi\`eme  partie du corollaire suivant, rappelons que si une vari\'et\'e ab\'elienne $A$ a r\'eduction semi-stable en une place $v$, la composante neutre de la fibre sp\'eciale est une extension d'une vari\'et\'e ab\'elienne par un tore; la dimension de ce tore s'appelle la {\it dimension torique}.   Remarquons \'egalement (cf. la preuve du th\'eor\`eme \ref{hallbis}) que si la vari\'et\'e ab\'elienne $A$ est de type I (resp. de type II) et poss\`ede mauvaise r\'eduction semi-stable en une place, alors la dimension torique de la fibre sp\'eciale est un multiple de $e$ (resp. de $2e$).
En rassemblant les cas o\`u l'on sait d\'emontrer qu'une vari\'et\'e ab\'elienne de type I ou II est pleinement de type Lefschetz (voir Th\'eor\`eme \ref{mtcases} et section \ref{aboutMT}) on obtient le corollaire suivant.

\begin{cor}\label{cormaintheo} Soit $E$ un corps de nombres totalement r\'eel de degr\'e $e=[E:\Q]$. Soit $A$ une vari\'et\'e ab\'elienne d\'efinie sur un corps de nombres, telle que $A$ est g\'eom\'etriquement simple de type I ou II et le centre de son alg\`ebre d'endomorphismes est $E$. On suppose de plus que l'une quelconque des trois hypoth\`eses suivantes est satisfaite :
\begin{enumerate}
\item La dimension relative $h$ est un nombre impair ou \'egal \`a 2.
\item On a $e=1$ (ou encore $E=\Q$) et la dimension relative $h$ n'appartient pas au sous-ensemble exceptionnel $\Sigma$.
\item La vari\'et\'e ab\'elienne $A$ est de type I (resp. II) et poss\`ede une place de mauvaise r\'eduction semi-stable avec dimension torique $e$ (resp. avec dimension torique $2e$).  
\end{enumerate}
On a alors
\begin{equation*}
\gamma(A)=\frac{2dhe}{1+2eh^2+he}=\frac{2\dim A}{\dim\MT(A)}.
\end{equation*}
\end{cor}

\medskip

\rem Lorsque $e=1$ et $A$ de type I, le point $2.$ correspond au cas ``g\'en\'erique" trait\'e dans \cite{HR}.

\medskip

\rem  Comme cela est d\'emontr\'e par Noot \cite{noot}, pour chaque type de groupe de Hodge, il existe des vari\'et\'es ab\'eliennes d\'efinies sur un corps de nombres, ayant ce groupe de Hodge donn\'e et v\'erifiant la conjecture de Mumford-Tate. En particulier, pour tout corps de nombres totalement r\'eel $E$ de degr\'e $e$ sur $\Q$ et tout entier $h$, il existe des vari\'et\'es ab\'eliennes $A$ d\'efinies sur un corps de nombres, de dimension $eh$ (resp. $2eh$) telles que $\End^0(A)=E$ (resp. $\End^0(A)$ est une alg\`ebre de quaternions ind\'efinie sur $E$) et qui sont pleinement de Lefschetz. Ces vari\'et\'es ab\'eliennes v\'erifient donc les hypoth\`eses du th\'eor\`eme.

\medskip

\rem
La preuve fournit une in\'egalit\'e l\'eg\`erement plus pr\'ecise que $|A(L)_{\rm tor}|\ll_{\epsilon}[L:K]^{\gamma(A)+\epsilon}$ de  la forme
\begin{equation*}
|A(L)_{\rm tor}|\ll [L:K]^{\gamma(A)+c_1/\log\log[L:K]}
\end{equation*}
On peut r\'ecrire ce r\'esultat en terme de minoration du degr\'e de l'extension engendr\'ee par un sous-groupe de torsion; nous donnons la teneur du r\'esultat ci-dessous dans le cas de l'extension engendr\'ee par un seul point de torsion.

\begin{theo}\label{unpoint} Soit $A/K$ une vari\'et\'e ab\'elienne  g\'eom\'etriquement simple de type I ou II, de dimension relative $h$ et pleinement de type Lefschetz. Il existe une constante $c_1:=c_1(A,K)>0$ telle que pour tout point de torsion $P$ d'ordre $n$ dans $A(\overline{K})$, on a:
\begin{equation*}
[K({P}):K]\geq c_1^{\omega(n)}n^{2h}.
\end{equation*}
\end{theo}

\rem On peut affaiblir l\'eg\`erement l'\'enonc\'e en \'ecrivant l'in\'egalit\'e sous la forme 
\[[K({P}):K]\gg n^{2h-\frac{c}{\log\log n}}\ \text{ ou encore }\ [K({P}):K]\gg_{\epsilon} n^{2h-\epsilon} \  (\forall\epsilon>0).\] 
\noindent Par ailleurs notons que, toujours pour un point $P$ d'ordre $n$, on a trivialement $[K({P}):K]\leq n^{2g}$. On voit donc que dans le cas totalement g\'en\'erique (appel\'e ``type GSp'' dans \cite{HR}) tous les points d'ordre $n$ engendrent des corps de degr\'es comparables  $n^{2g}\geq [K({P}):K]\gg n^{2g-\epsilon}$, mais qu'il n'en est plus de m\^eme quand $e\geq 2$.

\medskip 

\noindent Nous \'etendons les r\'esultats au cas   d'une vari\'et\'e ab\'elienne g\'eom\'etriquement isog\`ene \`a un produit  de vari\'et\'es ab\'eliennes de type I ou II.

\begin{theo}\label{theoproduit}
Soit $A/K$ une vari\'et\'e ab\'elienne isog\`ene sur $\bar K$ \`a un produit $A_1^{n_1}\times\dots\times A_d^{n_d}$ avec $A_i$ non $\bar K$-isog\`ene \`a $A_j$ pour $i\not= j$. On suppose que chaque facteur $A_i$ est de type I ou II et   est pleinement de type Lefschetz. Pour tout sous-ensemble non vide $I$ de $[1,d]$ on note
$A_I:=\prod_{i\in I}A_i$. On note $e_i$ la dimension du centre de $\End^0(A_i)$ et $h_i$ la dimension relative de $A_i$. Enfin on pose $d_i=1$ (resp. $d_i=2$) si $A_i$ est de type I (resp. de type II). On a alors
\begin{equation*}
\gamma(A)=\max_I\frac{2\sum_{i\in I}n_id_ih_ie_i}{1+\sum_{i\in I}2e_ih_i^2+h_ie_i}=\max_{I}\frac{2\sum_{i\in I}n_i\dim A_i}{\dim\MT(A_I)}.
\end{equation*}
\end{theo}

\begin{cor}\label{co11}
Soit $A$ une vari\'et\'e ab\'elienne g\'eom\'etriquement isog\`ene \`a un produit $A_1^{n_1}\times\dots\times A_d^{n_d}$ avec $A_i$ non $\bar K$-isog\`ene \`a $A_j$ pour $i\not= j$. On suppose que chaque facteur $A_i$ est de type I ou II et v\'erifie l'une quelconque des trois hypoth\`eses suivantes:
\begin{enumerate}
\item Le nombre  $h_i$ est  impair ou \'egal \`a 2.  
\item On a $e_i=1$ (ie $E_i=\Q)$ et $h_i\not\in\Sigma$.
\item La vari\'et\'e ab\'elienne $A_i$ est de type I (resp. II) et poss\`ede une place de mauvaise r\'eduction semi-stable avec dimension torique $e_i$ (resp. avec dimension torique $2e_i$).
\end{enumerate}
On a alors, avec les notations du th\'eor\`eme \ref{theoproduit}:
\begin{equation*}
\gamma(A)=\max_I\frac{2\sum_{i\in I}n_id_ih_ie_i}{1+\sum_{i\in I}2e_ih_i^2+h_ie_i}=\max_{I}\frac{2\sum_{i\in I}n_i\dim A_i}{\dim\MT(A_I)}.
\end{equation*}
\end{cor}

\rem Dans le contexte du corollaire \ref{co11}, l'analogue du th\'eor\`eme \ref{unpoint} dit simplement que:
\begin{equation*}
[K({P}):K]\geq c_1^{\omega(n)}n^{2h_0}
\end{equation*}
\noindent o\`u $h_0$ est le minimum des $h_i$.

\bigskip

\noindent{\bf R\'eductions}.  Le probl\`eme que nous \'etudions est clairement invariant par deux modifications : remplacer $K$ par une extension {\it finie} $K'$ et remplacer $A$ par une vari\'et\'e ab\'elienne $A'$ {\it isog\`ene} \`a $A$. Quitte \`a effectuer une extension finie de $K$ et \`a remplacer $A$ par une vari\'et\'e isog\`ene, nous pouvons donc supposer, et nous supposerons, dans le reste de l'article  les propri\'et\'es suivantes v\'erifi\'ees par $A/K$:
\begin{enumerate}
\item L'anneau des endomorphismes d\'efinis sur $K$ est \'egal \`a l'anneau des endomorphismes d\'efinis sur $\bar{K}$ ; on le notera donc $\End(A)$. 
\item L'anneau des endomorphismes $\End(A)$ est un ordre maximal dans $\End(A)\otimes\Q$.
\item La vari\'et\'e ab\'elienne s'\'ecrit $A=A_1^{n_1}\times\dots\times A_r^{n_r}$, avec $1\leq r$, $1\leq n_i$ et des vari\'et\'es ab\'eliennes $A_i$  absolument simples  non isog\`enes deux \`a deux.
\item L'adh\'erence de Zariski de $\rho_{\ell^{\infty},A}(G_K)$ est connexe.
\item Les repr\'esentations $\rho_{\ell^{\infty},A}$ sont ind\'ependantes sur $K$.
\end{enumerate}
En effet la possibilit\'e de l'obtention des trois premi\`eres propri\'et\'es par extension de $K$ et isog\'enie d\'ecoule des propri\'et\'es g\'en\'erales des vari\'et\'es ab\'eliennes, tandis que les points 4. et 5. s'obtiennent par une extension ad\'equate du corps $K$, en invoquant deux r\'esultats subtils de Serre  \cite{serrekr86,serrech}. Le point 5. 
est rappel\'e, ainsi que la d\'efinition de ``repr\'esentations ind\'ependantes" \`a la 
proposition \ref{indep}.    

\noindent Remarquons qu'on pourrait \'egalement songer \`a imposer que $A$ soit principalement polaris\'ee, mais cela forcerait \`a renoncer \`a la propri\'et\'e 2., il nous a sembl\'e plus commode de pr\'eserver cette derni\`ere.
 
\bigskip

\noindent{\bf Plan}. Le plan de ce texte est le suivant: dans la section suivante on rassemble un certain nombre de lemmes de th\'eorie des groupes et de combinatoire. Dans les sections 3 et  4 on d\'ecrit les accouplements $\lambda$-adiques d\'eduits de celui de Weil, ainsi que l'\'etude des propri\'et\'es des repr\'esentations galoisiennes qui sont utilis\'ees dans la suite. Les sections 6 et 7 contiennent les preuves des deux th\'eor\`emes cit\'es en introduction (th\'eor\`emes \ref{maintheo} et \ref{theoproduit}), d'abord dans le cas d'une vari\'et\'e ab\'elienne simple et pour un groupe annul\'e par $\ell$ puis dans le cas g\'en\'eral. Les preuves sont en fait \'ecrites pour $\ell$ assez grand et on indique dans la huiti\`eme section comment on peut modifier les preuves pour traiter les ``petits" nombres premiers de mani\`ere similaire. La courte neuvi\`eme section contient la d\'emonstration de la minoration du degr\'e de l'extension engendr\'ee par un point  de torsion (th\'eor\`eme \ref{unpoint}). La dixi\`eme et derni\`ere section est un appendice ind\'ependant du reste de l'article (sauf pour les notations), dans lequel nous montrons comment d\'eduire un \'enonc\'e du type ``conjecture forte de Mumford-Tate" de l'\'enonc\'e usuel de la conjecture de Mumford-Tate et expliquons comment prouver les quelques cas suppl\'ementaires de la conjecture de Mumford-Tate ne figurant pas dans la litt\'erature et \'enonc\'es dans les corollaires \ref{cormaintheo} et \ref{co11}.

\medskip

\noindent \textbf{Notations} Dans tout le reste de cet article, nous utiliserons les notations suivantes : 
\[ [L:K], \ (G:H),\ [G,G]\]
\noindent pour d\'esigner respectivement, le degr\'e de l'extension de corps $L$ sur $K$, l'indice du sous-groupe $H$ de $G$ dans $G$, le groupe des commutateurs de $G$. 

\medskip

\noindent Par ailleurs nous utiliserons les deux notations suppl\'ementaires suivantes : l'\textit{\'egalit\'e \`a indice fini pr\`es} $\que$ et la \textit{presque \'egalit\'e} $\pre$. Commen\c{c}ons par d\'efinir l'\'egalit\'e \`a indice fini pr\`es : si $L_1$, $L_2$ sont des corps de nombres contenus dans un corps $L$ (en pratique  nos corps seront tous contenus dans $K(A_{\rm tor})$ ou, si l'on pr\'ef\`ere dans $\bar{\Q}$ ou m\^eme $\C$) qui d\'ependent de $A/K$ et d'un autre ensemble de param\`etres $\Lambda$, nous \'ecrirons $L_1\que L_2$ pour dire qu'il existe une constante C(A/K) ne d\'ependant que de $A/K$ telle que les in\'egalit\'es $[L_1 : L_1\cap L_2] \leq C(A/K)$ et $[L_2 : L_1\cap L_2]\leq C(A/K)$ sont vraies pour toutes valeurs des param\`etres dans l'ensemble $\Lambda$. De m\^eme si $G_1$ et $G_2$ sont des sous-groupes d'un m\^eme groupe et d\'ependant d'un ensemble de param\`etres $\Lambda$, nous \'ecrirons $G_1\que G_2$ pour dire qu'il existe une constante C(A/K) ne d\'ependant que de $A/K$ telle que les in\'egalit\'es $(G_1 : G_1\cap G_2) \leq C(A/K)$ et $(G_2 : G_1\cap G_2)\leq C(A/K)$ sont vraies pour toutes valeurs des param\`etres dans l'ensemble $\Lambda$. Nous utiliserons la m\^eme notation pour des nombres. Ainsi,  si $N_1$, $N_2$ sont deux nombres (par exemple des cardinaux de groupes ou des  degr\'es d'extensions) $N_1\que N_2$ signifie qu'il existe deux constantes $C_1$ et $C_2$, ind\'ependantes des param\`etres, telles que $C_1N_1\leq N_2\leq C_2N_2$ (voir, par exemple, le lemme \ref{pmu1} pour une utilisation de cette derni\`ere notation).


\noindent Concernant la presque \'egalit\'e, dire que $E\pre F$ signifie par d\'efinition que  $E\que F$ et que de plus il y a en fait \'egalit\'e $E=F$ pour toutes les valeurs des param\`etres dans $\Lambda$ sauf \'eventuellement un nombre fini (d\'ependant \'eventuellement de $A/K$).

\medskip

Un exemple typique d'utilisation de la notation pr\'ec\'edente consiste par exemple (cf. proposition \ref{newprop}) \`a \'ecrire, pour $A/K$ une vari\'et\'e ab\'elienne simple de type I ou II de dimension relative $h$ et pleinement de type Lefschetz, que l'on a
\[ [\rho_{\ell}(G_K) , \rho_{\ell}(G_K)]\pre\Sp_{2h}(\O_E/\ell\O_E),\]
\noindent ceci signifiant qu'il existe une constante $\ell_0(A/K)$ d\'ependant de $A/K$ telle qu'il a \'egalit\'e pour tout $\ell\geq \ell_0(A/K)$ et que de plus pour tout $\ell$ les deux groupes sont commensurables, i.e.  l'indice de l'intersection des deux groupes dans l'un et l'autre est fini.

\medskip

\noindent Enfin, si $L_1$, $L_2$ sont des corps de nombres qui d\'ependent de $A/K$ et d'un autre ensemble de param\`etres $\Lambda$, nous \'ecrirons $[L_1:\Q]\ll [L_2:\Q]$ pour dire qu'il existe une constante C(A/K) ne d\'ependant que de $A/K$ telle que l'in\'egalit\'e $[L_1 : \Q] \leq C(A/K)[L_2:\Q]$ est vraie uniform\'ement sur $\Lambda$.

\medskip

\noindent\textbf{Remerciements.} Les auteurs remercient le rapporteur, dont la lecture attentive et les corrections et suggestions judicieuses ont permis d'am\'eliorer la pr\'esentation de ce texte.

\section{Lemmes de groupes\label{groupe}}

\noindent Soit $E/\Q$ un corps de nombres, d'anneau d'entiers $\O_E$. Si $\ell$ est un premier et $\lambda$ une place de $\O_E$ au dessus de $\ell$, on note $\O_{\lambda}$ le compl\'et\'e de $\O_E$ selon $\lambda$. De m\^eme on note $\F_{\lambda}$ le corps r\'esiduel correspondant.

\medskip

\noindent Nous rappelons maintenant des objets et notations provenant de \cite{HR} que nous utiliserons ensuite. Soit $V$ un $\Q$-espace vectoriel muni d'une forme symplectique (dans la suite du papier nous utiliserons ceci avec $V=H_1(A(\C),\Q)$). Dans les lemmes suivants nous introduisons $(e_1,\dots,e_{2g})$ une base symplectique de $V$. (\textit{i.e.} pour tout $1\leq i,j\leq g$, $e_i\cdot e_{g+i}=+1$ et $e_i\cdot e_j=0$ si $|i-j|\not= 0$).

\begin{lemme}\label{defprs} Soit $0\leq s\leq r\leq g$ avec $r\geq 1$. D\'efinissons $P_{r,s}$ le sous-groupe alg\'ebrique de $\Sp_{2g}$ fixant les vecteurs $e_1,\dots,e_r$ et les vecteurs  $e_{g+1},\dots,e_{g+s}$, c'est-\`a-dire
$$P_{r,s}:=\left\{M\in\Sp_{2g}\;|\; Me_i=e_i,\; i\in[1,r]\cup[g+1,g+s]\right\},$$
alors, $P_{r,s}$ est lisse sur $\O_E$ et sa codimension dans $\Sp_{2g}$ est :
\[\codim P_{r,s}=2sg+2rg-rs-\frac{r(r-1)}{2}-\frac{s(s-1)}{2}.\] 
\end{lemme}

\demo Il s'agit du lemme 2.24 de \cite{HR} (\'enonc\'e sur $\Z$ mais valable par la m\^eme preuve sur $\O_E$).$\findemo$

\begin{lemme}\label{lemmeA} Soit $\lambda$ une place de $\O_E$ au dessus d'un premier $\ell$ de $\Z$. En introduisant les groupes
\[ D_0:=\left\{\begin{pmatrix} I& 0\cr 0& \alpha I\cr\end{pmatrix}\in\GL_{2g}(\O_{\lambda})\ |\ \alpha\in \Z_{\ell}^{\times}\right\},\text{ et } G_{\lambda}:=\left\{M\in\GSp_{2g}(\O_{\lambda})\ |\ \mult(M)\in \Z_{\ell}^{\times}\right\},\] 
\noindent on a
\[D_0\cdot\Sp_{2g}(\O_{\lambda})=G_{\lambda}.\]
\noindent Le m\^eme \'enonc\'e vaut en rempla\c{c}ant $\O_{\lambda}$ par $\F_{\lambda}$ et $\Z_{\ell}$ par $\F_{\ell}$.
\end{lemme}
\demo  Il s'agit du lemme 2.12 de \cite{HR} dans sa version $\lambda$-adique : soit $M\in G_{\lambda}$ de multiplicateur $\mult(M)$. La matrice
$\begin{pmatrix} I_g& 0\cr 0& \lambda(M)^{-1} I_g\cr\end{pmatrix}M$ est dans $\Sp_{2g}(\O_{\lambda})$. $\findemo$

\medskip

\begin{lemme}\label{cle}  Soit $G_E$ un sous-groupe alg\'ebrique sur $E$ de $\GL_E$. Soit $t\in\N^*$ et soit $\mathcal{G}_{1,E},\ldots,\mathcal{G}_{t,E}$ une suite de sous-groupes alg\'ebriques de $G_E$. On note $G$ (respectivement $\mathcal{G}_{i}$) l'adh\'erence de Zariski de $G_E$ (respectivement de $\mathcal{G}_{i,E}$) dans $\GL_{\O_E}$  sur $\O_E$. Il existe des constantes $C_1>0$, $C_2>0$ telles que la propri\'et\'e suivante est vraie : soient $\ell$ un premier de $\Z$ et $\lambda$ une place de $\O_E$ au-dessus de $\ell$. Soient $G_1\subset G_{2}\subset\dots\subset G_t$ une suite de sous-groupes alg\'ebriques sur $\O_{\lambda}$ de $G_{\O_{\lambda}}$. On suppose que pour tout $i$, le groupe $G_i$ est conjugu\'e sur $\F_{\lambda}$ \`a $\mathcal{G}_i$. On note $g_i:=\dim\mathcal{G}_i=\dim G_i$ et $d_i:=\codim_G \mathcal{G}_i=\codim_{G_{\O_{\lambda}}}G_i$ et on pose, pour toute suite croissante d'entiers   
$0=m_0<m_1<m_2<\dots<m_t$~:
\[H(m_1,\dots,m_t)=\left\{M\in G(\O_{\lambda})\;|\; M\in G_i\mod \lambda^{m_i}\right\}.\]
\noindent Pour tous les $\ell$ tels que $G$ et les $\mathcal{G}_i$ sont lisses sur $\F_{\lambda}$, on a alors
\[C_1\times\left(G(\O_{\lambda}): H(m_1,\dots,m_t)\right)\geq\text{Card}(\F_{\lambda})^{\sum_{i=1}^td_i(m_i-m_{i-1})}\geq C_2\times\left(G(\O_{\lambda}): H(m_1,\dots,m_t)\right).\]
\end{lemme}
\demo Il s'agit du lemme 2.4 de \cite{HR} dont la preuve reste valable, en rempla\c{c}ant $\Z$ par $\O_E$, ainsi que $\ell$ par $\lambda$ et $\F_{\ell}$ par $\F_{\lambda}$.$\findemo$

\medskip

\noindent Nous donnons dans ce qui suit l'analogue d'un lemme prouv\'e pour $\SL_2(\Z_{\ell})$ par Serre. L'\'enonc\'e en vue est le suivant, o\`u les notations suivantes sont utilis\'ees.

\medskip

\defi Nous dirons qu'une sous-alg\`ebre de Lie de $gl_m$ \textit{poss\`ede la propri\'et\'e $(\mathcal{CN})$ (des carr\'es nuls)} si elle est engendr\'ee, comme espace vectoriel, par des matrices de carr\'e nul.

\medskip

\exemple Les alg\`ebres $\mathfrak{sl}_{m}$ et $\mathfrak{sp}_{2m}$ ont la propri\'et\'e $(\mathcal{CN})$ mais pas $\mathfrak{so}_m$. Les matrices $E_{ij}$ ayant un seul coefficient non nul (sur la $i$-\`eme ligne et $j$-\`eme colonne) sont de carr\'e nul. En dimension 2, on a  de m\^eme que $\begin{pmatrix}a&1\cr -a^2&-a\cr\end{pmatrix}$ est de carr\'e nul, donc $$\begin{pmatrix}a&b\cr c&-a\cr\end{pmatrix}=\begin{pmatrix}a&1\cr -a^2&-a\cr\end{pmatrix}+\begin{pmatrix}0&b-1\cr 0&0\cr\end{pmatrix}+\begin{pmatrix}0&0\cr c+a^2&0\cr\end{pmatrix}$$
 est bien somme de matrices de carr\'e nul. Quand $m$ est quelconque, on en tire ais\'ement que les matrices diagonales de trace nulle sont sommes de matrices de carr\'e nul. L'alg\`ebre $\mathfrak{sp}_{2m}$ est l'alg\`ebre des matrice $\begin{pmatrix} A& B\cr C& D\cr\end{pmatrix}$ telles que $D+{}^tA=0$, $B$ et $C$ sont sym\'etriques. On \'ecrit ais\'ement une telle matrice comme somme de matrices:
 $$\begin{pmatrix} 0& S\cr 0& 0\cr\end{pmatrix},\quad \begin{pmatrix} 0& 0\cr S& 0\cr\end{pmatrix}\quad \begin{pmatrix} \delta& \delta\cr -\delta& -\delta\cr\end{pmatrix},\quad
 \begin{pmatrix} U& 0\cr 0& -^tU\cr\end{pmatrix}$$
 o\`u $S$ est sym\'etrique, $U$ a trace nulle et $\delta$ est la matrice dont le seul terme non nul est dans le coin sup\'erieur gauche et vaut 1; les trois premi\`eres sont de carr\'e nul et la quatri\`eme est somme de matrices de carr\'es nuls d'apr\`es la propri\'et\'e pour $\mathfrak{sl}_m$.
Enfin notons qu'une matrice anti-sym\'etrique de carr\'e nul est elle-m\^eme nulle, donc $\mathfrak{so}_n$ ne poss\`ede pas la propri\'et\'e des carr\'es nuls.

\medskip

\noindent{\bf Notations.} On note   $\mathcal{O}$ l'anneau d'entiers d'un corps $p$-adique, $\varpi$ une uniformisante, $\F=\mathcal{O}/\varpi\mathcal{O}$ le corps r\'esiduel et $e$ l'indice de ramification, i.e. $p =\varpi^eu$ avec $u\in\mathcal{O}^{\times}$.

\medskip

\begin{lemme}\label{releveserre}  Soit $G$ un sous-groupe alg\'ebrique lisse de $\GL_m/\mathcal{O}$ et $H$ un sous-groupe ferm\'e de $G(\mathcal{O})$. Consid\'erons, pour tout entier $n\geq 1$, les applications $\pi_n:H\rightarrow\GL_m(\mathcal{O}/\varpi^n\mathcal{O}) $, alors on a
\begin{enumerate}
\item Si $\pi_{e+1}(H)=G(\mathcal{O}/\varpi^{e+1}\mathcal{O})$ et $p\geq e+2$ alors $H=G(\mathcal{O})$; si $p\leq e+1$ et $\pi_{m}$ surjective avec $m\geq \frac{ep+1}{p-1}$,  alors $H=G(\mathcal{O})$;
\item Si $p\geq 2e+3$, $\pi_e(H)=G(\mathcal{O}/\varpi)$ et $Lie(G_{\F})$ a la propri\'et\'e $(\mathcal{CN})$ alors $H=G(\mathcal{O})$;
\end{enumerate}
\noindent En particulier, lorsque $G=\SL_m$ ou $\Sp_m$, si $p\geq 5$ et $K/\Q_p$ non ramifi\'e, alors $\pi_1(H)=G(\F)$ entra\^{\i}ne $H=G(\mathcal{O})$ .
\end{lemme}
\demo Notons $\mathcal{L}:=Lie(G_{\F})$. Commen\c{c}ons par observer que, d'apr\`es l'hypoth\`ese de lissit\'e, si une matrice s'\'ecrit $M=I+\varpi^nB$ alors $\pi_{n+1}(M)\in G(\O/\varpi^{n+1}\O)$ si et seulement si $\pi_1(B)\in \mathcal{L}$. On prouve maintenant par r\'ecurrence sur $n$ que la projection $H\rightarrow G(\O/\varpi^n\O)$ est surjective et donc que $H$ est dense dans 
$G(\O)$ et donc \'egal \`a ce dernier. Supposons la propri\'et\'e vraie au cran $n$ et montrons-la au cran $n+1$. Soit donc $A\in G(\O)$, on sait donc qu'il existe $A_1\in H$ telle que $A\equiv A_1[\varpi^n]$ et, quitte \`a remplacer $A$ par $AA_1^{-1}$ on peut supposer que $A\equiv I[\varpi^n]$ soit encore $A=I+\varpi^nB$ avec donc $\pi_1(B)\in \mathcal{L}$. 
 Par hypoth\`ese, il existe $Z\in H$ telle que $Z\equiv I+\varpi^{n-e}B[\varpi^{n-e+1}]$
ou encore $Z=I+\varpi^{n-e}B+\varpi^{n-e+1}C$. Posons $Y:=Z^{p}$ alors
$$Y=I+p\varpi^{n-e}(B+\varpi C)+\sum_{h=2}^{p-1}{p\choose h}\varpi^{h(n-e)}(B+\varpi C)^h +\varpi^{p(n-e)}(B+\varpi C)^{p}$$
On a $e+h(n-e)\geq n+1$ si $n\geq e+1$ et $p(n-e)\geq n+1$ si $n\geq \frac{ep+1}{p-1}$; ainsi si $p\geq e+2$ on voit que 
$Y\equiv I+\varpi^{n}uB\mod \varpi^{n+1}$. On peut bien s\^ur refaire ce calcul en rempla\c{c}ant $B$ par $u^{-1}B$ et conclure. La surjectivit\'e de $\pi_{e+1}$ suffit donc pour entra\^{\i}ner $H=G(\O)$ si $p\geq e+2$ (resp. la surjectivit\'e de $\pi_{m}$ avec $m\geq \frac{ep+1}{p-1}$ pour $p\leq e$).

Si maintenant $n=e$, on reprend le calcul en supposant d'abord que $\pi_1(B^2)=0\in\mathcal{L}$. Observons  que, dans ce cas,  $B^{2j}\equiv 0\mod \varpi^j$ donc 
$$(B+\varpi C)^p=B^p+\varpi\left(B^{p-1}C+\dots\right)+\dots+\varpi^r\left(\dots +B^{j_0}CB^{j_1}C\dots CB^{j_s}\right)+\dots+\varpi^pC^p$$
 avec $B^p\equiv 0\mod \varpi^{(p-1)/2}$ et $\varpi^rB^{j_0}CB^{j_1}C\dots CB^{j_s}\equiv 0\mod \varpi^m$ avec
 $$m\geq r+\left[\frac{j_0}{2}\right]+\dots \left[\frac{j_s}{2}\right]\geq r+\sum_{i=0}^s(\frac{j_i}{2}-\frac{1}{2})=r+\frac{p-r}{2}-\frac{s+1}{2}\geq\frac{p-1}{2}$$
 
Si maintenant $Z=I+B+\varpi C$ et $Y=Z^p$ alors
$$Y=I+p(B+\varpi C)+\sum_{h=2}^{p-1}{p\choose h}\varpi^h(B+\varpi C)^{p-h}+(B+\varpi C)^p\equiv I+pB+(B+\varpi C)^p\mod\varpi^{e+1}$$
La condition $(p-1)/2\geq e+1$ \'equivaut \`a $p\geq 2e+3$.
On trouve ainsi un \'el\'ement $Y\in H$ tel que $Y=I+\varpi^eB\mod\varpi^{e+1}$.  Pour le cas g\'en\'eral, on aura $B=B_1+\dots+B_s$ avec 
$\pi_1(B_i^2)=0\in\mathcal{L}$ et $I+\varpi^e B\equiv(I+\varpi^eB_1)\dots(I+\varpi^e B_s)\mod \varpi^{e+1}$, ce qui permet de conclure. $\findemo$

\medskip

\noindent Le lemme \'el\'ementaire suivant est d\'emontr\'e dans \cite{HR} (lemme 2.8).

\begin{lemme}\label{combielem2} Soit $d\geq 1$ un entier, et pour tout $i\in\{1,\ldots,d\}$, soient $t_i\geq 1$ des entiers. Pour $i\leq d$ et $j\leq t_i$, on se donne \'egalement des entiers $a_{ij}$ et $b_{ij}$, strictement positifs. On a l'\'egalit\'e
\begin{equation}\label{formulelem2}
\sup_{\underset{1\leq i\leq d}{m_{i1}\geq\dots\geq m_{i t_i}}}\left\{\frac{\sum_{i=1}^d\sum_{j=1}^{t_i}a_{ij}m_{ij}}{\sum_{i=1}^d\sum_{j=1}^{t_i}b_{ij}m_{ij}}\right\}=
\max_{\underset{1\leq i\leq d}{1\leq h_i\leq t_i}}\left\{\frac{\sum_{i=1}^d\sum_{j=1}^{h_i}a_{ij}}{\sum_{i=1}^d\sum_{j=1}^{h_i}b_{ij}}\right\},
\end{equation}
\noindent le sup dans le membre de gauche \'etant pris sur les entiers $m_{ij}$ ordonn\'es pour $1\leq i\leq d$ par $m_{i1}\geq\dots\geq m_{i t_i}$ et tels que $m_{i1}\not=0$.
\end{lemme}

\section{Repr\'esentations et accouplements $\lambda$-adiques}

\noindent L'\'etude de la repr\'esentation galoisienne ad\'elique se ram\`ene essentiellement \`a l'\'etude des repr\'esen\-ta\-tions $\ell$-adiques, gr\^ace au r\'esultat suivant d\^u \`a Serre.

\medskip

\defi Une famille de repr\'esentations $(\rho_i)_{i\in I}:G_K\rightarrow \prod_{i\in I}\GL(V_i)$ index\'ee par un ensemble $I$, est dite \textit{ind\'ependante} si 
\[(\rho_i)_{i\in I}(G_K)=\prod_{i\in I}\left(\rho_i(G_K)\right).\]
\noindent Les $(\rho_i)_{i\in I}$ sont \textit{presque ind\'ependantes} s'il existe une extension finie $K'/K$ telles que les restrictions \`a $G_{K'}$ sont ind\'ependantes.

\begin{prop}\label{indep} \textnormal{\textbf{(Serre)}}
Les repr\'esentations $\ell$-adiques associ\'ees \`a une vari\'et\'e ab\'elienne $A$ sur un corps de nombres $K$ sont presque ind\'ependantes.
\end{prop}
\demo C'est le th\'eor\`eme 1 de \cite{college8586}, cf. \'egalement \cite{serrech} th\'eor\`eme 1 et paragraphe 3.1.\hfill $\Box$

\medskip

\noindent Cet \'enonc\'e se traduit concr\`etement en disant que, quitte \`a remplacer $K$ par une extension finie, pour tous premiers $\ell_1$, $\ell_2$ distincts, les extensions $K(A[\ell_1^{\infty}])/K$ et $K(A[\ell_2^{\infty}])/K$ sont lin\'eairement disjointes.

\medskip

\noindent Nous utiliserons en parall\`ele les d\'ecompositions $\ell$-adiques et $\lambda$-adiques correspondant aux types I et II; ces d\'ecompositions s'\'ecrivent pour tout $\ell$ au niveau des $\Q_{\ell}$-repr\'esentations $V_{\ell}(A)$ et pour $\ell$ assez grand (hors d'un ensemble fini de $\ell$) pour les $\Z_{\ell}$-repr\'esentations $T_{\ell}(A)$.

\subsection{Modules de Tate, repr\'esentations $\ell$-adiques et $\lambda$-adiques}
\noindent Soit $\ell$ un premier quelconque. On consid\`ere dans toute la suite de cette section 3. une vari\'et\'e ab\'elienne $A/K$ g\'eom\'etriquement simple de type I ou II, telle que $\End_K(A)=\End_{\bar K}(A)$ et telle que $\End_K(A)$ est un ordre maximal de $D:=\End_K(A)\otimes_{\Z}\Q$. Nous noterons $E$ le centre de $D$.

\medskip

\noindent On a la d\'ecomposition suivante : $\O_{\ell}=\prod_{\lambda |\ell}\O_{\lambda}$ o\`u $\O_{\lambda}$ est l'anneau des entiers du compl\'et\'e $E_{\lambda}$ de $E$ pour la place $\lambda$. En notant $f(\lambda):=[E_{\lambda}:\Q_{\ell}]$ et $e(\lambda)$ le degr\'e de ramification de $\lambda |\ell$, on a $\sum_{\lambda|\ell}e(\lambda)f(\lambda)=e=[E:\Q]$.

\medskip

\noindent  Le module de Tate $\ell$-adique, $T_{\ell}(A)$, est muni d'une action de $\O_{\ell}$ et se d\'ecompose en 

\[T_{\ell}(A)=\prod_{\lambda|\ell}\T_{\lambda} \ \text{ o\`u }\ \T_{\lambda}:=T_{\ell}(A)\otimes_{\O_{\ell}}\O_{\lambda}.\] 

\noindent En inversant $\ell$ on obtient les espaces $V_{\ell}(A)$ et $V_{\lambda}$ \`a partir de $T_{\ell}(A)$ et de $\T_{\lambda}$ : 
\[V_{\ell}(A):=T_{\ell}(A)\otimes_{\Z_{\ell}}\Q_{\ell}\ \text{ et }\ V_{\lambda}:=\T_{\lambda}\otimes_{\O_{\lambda}}E_{\lambda}.\]

\noindent La repr\'esentation $\ell$-adique \'etant $\O_E$-lin\'eaire, elle se d\'ecompose diagonalement selon les $\lambda$ divisant $\ell$ en repr\'esentations $\lambda$-adiques (cf. \cite{ribetrm} paragraphe II) :
\[\rho_{\ell^{\infty}}=(\prod_{\lambda |\ell}\rho_{\lambda^{\infty}}) : G_K\rightarrow \prod_{\lambda |\ell }\Aut(\T_{\lambda}).\]

\noindent La repr\'esentation $\ell$-adique modulo $\ell$ \'etant $\O_E/\ell\O_E$-lin\'eaire, elle se d\'ecompose par r\'eduction modulo $\ell$ diagonalement selon les $\lambda$ en repr\'esentations $\lambda$-adiques :
\[\rho_{\ell}=(\prod_{\lambda |\ell}\rho_{\lambda}) : G_K\rightarrow \Aut_{\O_{\ell}/\ell\O_{\ell}}(A[\ell])=\prod_{\lambda |\ell }\Aut_{\O_{\lambda}/\ell\O_{\lambda}}(\T_{\lambda}/\ell\T_{\lambda}).\]
\noindent 
\noindent Nous noterons dans toute la suite $G_{\lambda}$ le groupe de Galois correspondant \`a $\rho_{\lambda}$ (il est a priori \`a valeurs dans $\GL_{2h}(\O_{\lambda}/\ell\O_{\lambda})$). La m\^eme chose vaut au niveau $\ell$-adique et on sait (cf. par exemple \cite{Chi} p.319) que ces repr\'esentations $\lambda$-adiques sont munies naturellement d'un accouplement de Weil $\lambda$-adique provenant de l'accouplement $\ell$-adique. Nous rappelons ceci dans le paragraphe suivant. 

\subsection{Accouplements $\ell$-adique et $\lambda$-adique}
\noindent Nous supposons ici de plus que $A$ est polaris\'ee par une polarisation $\phi$ et que $\ell$ est un premier ne divisant pas $\deg(\phi)$.

\medskip

\noindent Rappelons la construction de l'accouplement $\lambda$-adique (cf. par exemple \cite{BGK} paragraphes 3 et 4). On commence pour cela par l'accouplement de Weil $\ell$-adique usuel : 
\[\phi_{\ell^{\infty}} : T_{\ell}(A)\times T_{\ell}(A)\rightarrow \Z_{\ell}(1)=\varprojlim \mu_{\ell^m}.\]
\noindent L'accouplement usuel de Weil est non-d\'eg\'en\'er\'e (modulo $\ell^n$ pour tout $n\geq 1$) car $\ell$ ne divise pas $\deg(\phi)$. De plus, si $\dagger$ d\'esigne l'involution de Rosati sur $\End^0(A)$ associ\'ee \`a la polarisation d\'efinissant l'accouplement on aura 
$\phi_{\ell^{\infty}}(ax,y)=\phi_{\ell^{\infty}}(x,a^{\dagger}y)$ pour $x,y\in T_{\ell}(A)$ et $a\in \O_{\ell}$.

\begin{lemme}\label{pairing}Notons $\O^{\star}_{\ell}$ le dual de $\O_{\ell}$ pour la dualit\'e donn\'ee par la trace $Tr_{E_{\ell}/\Q_{\ell}}$. Il existe un unique accouplement $\O_{\ell}$-lin\'eaire, $\phi_{\ell^{\infty}}^{\star}: T_{\ell}(A)\times T_{\ell}(A)\rightarrow \O^{\star}_{\ell}(1),$ tel que 
\[Tr_{E_{\ell}/\Q_{\ell}}(\phi_{\ell^{\infty}}^{\star})=\phi_{\ell^{\infty}}.\]
\end{lemme}
\demo Il s'agit essentiellement du lemme 3.1 de \cite{BGK} (cf. aussi sublemma 4.7 de \cite{deligne}, page 55). Le preuve est la suivante : il s'agit de v\'erifier que le morphisme donn\'e par la trace, de $\Hom_{\O_{\ell}}\left(T_{\ell}(A)\otimes_{\O_{\ell}}T_{\ell}(A),\O_{\ell}^{\star}\right)$ vers $\Hom_{\Z_{\ell}}\left(T_{\ell}(A)\otimes_{\O_{\ell}}T_{\ell}(A),\Z_{\ell}\right)$ est un isomorphisme. Or ces deux objets sont des $\Z_{\ell}$-modules libres de m\^eme rang et la preuve du lemme 3.1 de \cite{BGK} donne la surjectivit\'e. \hfill $\Box$

\medskip

\noindent \textbf{Hypoth\`ese : } On suppose dans la fin de ce paragraphe que $\ell$ est de plus non ramifi\'e dans $E/\Q$, ie que $\O_{\ell}^{\star}=\O_{\ell}$.

\medskip

\noindent On a dans ce cas l'accouplement $\O_{\ell}$-lin\'eaire $\phi_{\ell^{\infty}}^{\star}: T_{\ell}(A)\times T_{\ell}(A)\rightarrow \O_{\ell}(1)$. Par projection on construit alors l'accouplement $\lambda$-adique, $\O_{\lambda}$-lin\'eaire, de la fa\c{c}on suivante :
\[\phi_{\lambda^{\infty}} : \T_{\lambda}\times \T_{\lambda}\rightarrow \O_{\lambda}(1)\]
\noindent tel que $\phi_{\ell^{\infty}}^{\star}(x,y)\otimes 1=\phi_{\lambda^{\infty}}(x\otimes 1,y\otimes 1).$

\medskip

\noindent Tout comme l'accouplement $\ell$-adique, l'accouplement $\lambda$-adique est Galois \'equivariant : 
\[\forall \sigma\in G_K,\forall x,y\in \T_{\lambda},\ \ \phi_{\lambda^{\infty}}(\sigma x,\sigma y)=\phi_{\lambda^{\infty}}(x,y)^{\sigma},\]
\noindent l'action de Galois \`a gauche se faisant via la repr\'esentation $\lambda$-adique, et \`a droite, via le caract\`ere cyclotomique $\ell$-adique usuel.

\medskip

\noindent Nous noterons enfin $\phi_{\ell^{\infty}}^0$ et $\phi_{\lambda^{\infty}}^0$ les accouplements $\ell$-adiques et $\lambda$-adiques, d\'efinis de mani\`ere similaires sur $V_{\ell}(A)$ et $V_{\lambda}$, \`a valeurs respectivement dans $E_{\ell}$ et $E_{\lambda}$. Ils sont d\'efinis sans restriction pour tout $\ell$.

\subsection{Galois pour les vari\'et\'es de type I et II}
\noindent Soit $\ell$ un premier quelconque.

\medskip

\noindent Dans le cas de type I, les $(V_{\lambda},\phi_{\lambda^{\infty}}^0)$ fournissent des repr\'esentations irr\'eductibles symplectiques. Dans le cas de type II, on a une d\'ecomposition plus fine (cf. \cite{Chi2}, \cite{milne99}, \cite{BGK})~:
\begin{equation*}
V_{\lambda}=W_{\lambda}(A)\oplus W_{\lambda}(A),
\end{equation*}
o\`u $(W_{\lambda}(A),\phi^0_{\lambda^{\infty}|W_{\lambda}(A)})$ est maintenant irr\'eductible symplectique. Cependant, comme nous le d\'etaillons plus loin, cette d\'ecomposition,  dans le cas de type II, n\'ecessite, pour le nombre fini de  premiers $\ell$ ``ramifi\'es" pour l'alg\`ebre de quaternions,  d'\'etendre les scalaires \`a une extension quadratique.

\medskip

\noindent Nous avons $ V_{\lambda}:=W_{\lambda}(A)\oplus W_{\lambda}(A)$ si $A$ est de type II et nous posons $W_{\lambda}(A):=V_{\lambda}$ si $A$ est de type I. Autrement dit en notant $d=1$ si $A$ est de type I et $d=2$ si $A$ est de type II, nous aurons toujours 
\[V_{\lambda}=W_{\lambda}(A)^d,\]
\noindent et dans tous les cas, la repr\'esentation $W_{\lambda}(A)$ est irr\'eductible, symplectique et, lorsque $A$ est de type I  (resp. de type II), le module $V_{\ell}(A)$ est isomorphe \`a
$\prod_{\lambda}W_{\lambda}(A)$ (resp. \`a la somme de deux copies de ce produit). 
\noindent On pose ensuite dans tout les cas
\[T_{\lambda}(A):=\T_{\lambda}\cap W_{\lambda}(A).\]

\medskip

\defi \label{ex} Nous noterons dans la suite $S_{\text{ex}}(A)$  l'ensemble fini des $\ell$ divisant le degr\'e de la polarisation fix\'ee $\phi$ de $A$, des $\ell$ ramifi\'es dans $\O_E$ et, dans le cas de type II, des $\ell$ tels que l'alg\`ebre de quaternions $D$ est non d\'ecompos\'ee en au moins un $\lambda|\ell$. 

\medskip

\noindent Au niveau des $\Z_{\ell}$-modules, la d\'ecomposition perdure, au moins pour $\ell$ hors de $S_{\text{ex}}(A)$. Nous rassemblons ces \'enonc\'es dans la proposition suivante.

\begin{prop}\label{2copies} Soit $W_{\lambda}(A)$ le $E_{\lambda}$-module galoisien symplectique d\'efini ci-dessus, qu'on identifie \`a un sous-module de $V_{\ell}(A)$.
\begin{enumerate}
\item La repr\'esentation $W_{\lambda}(A)$ est irr\'eductible et symplectique.
\item Si $A$ est de type I (resp. de type II), on a une d\'ecomposition $V_{\ell}(A)=\prod_{\lambda|\ell}W_{\lambda}(A)$, resp. $V_{\ell}(A)=\prod_{\lambda|\ell}\left(W_{\lambda}(A)\oplus W_{\lambda}(A)\right)$. Toutefois, dans le cas de type II et pour les premiers ramifi\'es de l'alg\`ebre de quaternions, cette d\'ecomposition ne s'obtient qu'apr\`es tensorisation par   une extension quadratique de $E$.
\item Pour $\ell\notin S_{ex}(A)$, on a une d\'ecomposition analogue pour les $\mathcal{O}\otimes\Z_{\ell}$-modules : 
$T_{\ell}(A)=\prod_{\lambda|\ell}T_{\lambda}(A)$ si $A$ est de type I, resp. $T_{\ell}(A)=\prod_{\lambda}\left(T_{\lambda}(A)\oplus T_{\lambda}(A)\right)$ si $A$ est de type II.
\end{enumerate}
\end{prop}

\demo Pour $A$ de  type I, voir \cite{ribetrm}. Pour $A$ de type II, cet \'enonc\'e est prouv\'e pour $\ell$ assez grand sur $\Q_{\ell}$ dans \cite{Chi2} et sur $\Z_{\ell}$ dans \cite{BGK}, le point essentiel  \'etant l'identification de $\End^0(A)\otimes\Q_{\ell}$ \`a un produit d'alg\`ebres de matrices $\prod_{\lambda}M_2(E_{\lambda})$, ce qui est possible justement quand $\ell$ n'est pas ramifi\'e. L'alg\`ebre de quaternions est toujours d\'eploy\'ee sur une extension quadratique et l'on peut donc se ramener au cas pr\'ec\'edent apr\`es tensorisation par une telle extension; plus pr\'ecis\'ement, en choisissant $F$ extension quadratique de $E$ telle que
$D\otimes_E F\cong M_2(F)$ et en notant $F_{\lambda}:=E_{\lambda}\otimes_E F$ et $V_{\ell}\otimes_EF=\prod_{\lambda}V_{\lambda,F}$, chaque $V_{\lambda,F}$ est un module sur $D\otimes_EF_{\lambda}\cong M_2(F_{\lambda})$ et on peut alors imiter le proc\'ed\'e d\'ecrit dans {\it loc. cit.}.
Remarquons que la preuve sur $\Q_{\ell}$ pour tout $\ell$ peut  aussi \^etre extraite de \cite{milne99}. 

\medskip

\defi D\'efinissons la {\it dimension relative} de $A$ simple avec $D:=\End^0(A)$ et $E$ le centre de $D$, par la formule:
\begin{equation*}
\dimrel(A):=\frac{\dim A}{[E:\Q]\sqrt{[D:E]}}
\end{equation*}
\noindent Ainsi, si $e=[E:\Q]$,  la dimension relative $h$ d'une vari\'et\'e ab\'elienne de dimension $g$ de type I (resp. de type II) est $h=\frac{g}{e}$ (resp. $h=\frac{g}{2e}$).

\medskip

\noindent On a alors dans le cas I et II une repr\'esentation irr\'eductible symplectique de dimension $2h$, et une inclusion
\[H_{\ell,A}\subset \prod_{\lambda|\ell}\Sp(W_{\lambda}(A),\phi_{\lambda^{\infty}}^0)\cong \prod_{\lambda|\ell}\Sp_{2h,E_{\lambda}},\]
\noindent qu'il convient de mettre en parall\`ele avec l'inclusion
\[\Hdg(A)_{\Q_{\ell}}\subset(\Res_{E/\Q}\Sp_{2h,E})_{\Q_{\ell}}.\]

\noindent Nous allons nous placer dans le cas g\'en\'erique o\`u l'on a \'egalit\'e dans les deux inclusions pr\'ec\'edentes (la premi\`ere impliquant d'ailleurs la seconde
puisque $H_{\ell,A}\subset \Hdg(A)_{\Q_{\ell}}$). Le th\'eor\`eme suivant pr\'ecise des conditions o\`u l'on sait que l'\'egalit\'e voulue est toujours vraie.

\begin{theo}\label{mtcases} \textnormal{\textbf{(Banaszak-Gajda-Kraso\'n \cite{BGK}, Pink \cite{pink}, Hall \cite{hall})}} \label{mt12} Soit $A/K$ une vari\'et\'e ab\'e\-lienne de type I ou II, de dimension relative $h$; notons $E$ le centre de $\End^0(A)$. Supposons de plus l'une des trois conditions suivantes r\'ealis\'ee:
\begin{enumerate}
\item  L'entier  $h$ est impair ou \'egal \`a 2,
\item On a $E=\Q$ et $h$ n'appartient pas \`a l'ensemble exceptionnel $\Sigma$ d\'efini en (\ref{defdeS}),
\item La vari\'et\'e ab\'elienne $A$ est de type I (resp. II) et poss\`ede une place de mauvaise r\'eduction semi-stable avec dimension torique $e$ (resp. avec dimension torique $2e$).
\end{enumerate}
on a alors $H_{\ell,A}=\prod_{\lambda}\Sp_{2h,E_{\lambda}}.$
\end{theo}
\demo Le r\'esultat est d\'emontr\'e dans \cite{BGK} sous l'hypoth\`ese de l'alin\'ea 1 (cf. \cite{lombardo} remark 2.25 pour le cas $h=2$) ; le r\'esultat est d\'emontr\'e explicitement dans \cite{pink} sous l'hypoth\`ese de l'alin\'ea 2, pour une vari\'et\'e ab\'elienne de type I, mais on peut extraire de \cite{pink} une preuve pour le type II ; nous indiquons comment en appendice de cet article (Cf section \ref{aboutMT}). Le r\'esultat est d\'emontr\'e  dans \cite{hall} sous l'hypoth\`ese de l'alin\'ea 3, pour une vari\'et\'e ab\'elienne de type I v\'erifiant $e=1$ ; nous indiquons  en appendice de cet article comment \'etendre les arguments de \cite{hall} aux cas \'enonc\'es.\hfill $\Box$

\subsection{Variantes $\lambda$-adiques modulo $\lambda^n$}
\noindent On suppose ici que $\ell$ est un premier quelconque. Soit $n\geq 1$ un entier. On a par r\'eduction modulo $\ell^n$,
\[A[\ell^n]=T_{\ell}(A)\otimes_{\Z_{\ell}}\Z/\ell^n\Z=T_{\ell}(A)/\ell^nT_{\ell}(A).\]
\noindent De m\^eme par r\'eduction modulo $\lambda^n$, on voit que 
\[A[\lambda^n]:=\T_{\lambda}\otimes_{\O_{\lambda}}\O_{\lambda}/\lambda^n=\T_{\lambda}/\lambda^n\T_{\lambda}.\]
\noindent On pose ensuite
\[T_{\lambda}[\lambda^n]:=T_{\lambda}(A)\otimes_{\O_{\lambda}}\O_{\lambda}/\lambda^n=T_{\lambda}(A)/\lambda^nT_{\lambda}(A).\]
\noindent En utilisant que $\ell^n\O_{\lambda}=\lambda^{e(\lambda)n}$ on voit que par r\'eduction modulo $\ell^n$ on obtient :
\[A[\ell^n]=T_{\ell}(A)/\ell^nT_{\ell}(A)=\prod_{\lambda|\ell}\T_{\lambda}/\ell^n\T_{\lambda},\] 
\noindent et
\[T_{\lambda}(A)/\ell^nT_{\lambda}(A)=T_{\lambda}(A)\otimes_{\O_{\lambda}}\O_{\lambda}/\lambda^{e(\lambda)n}=T_{\lambda}[\lambda^{e(\lambda)n}].\]

\noindent Soit $\pi_{\lambda}$ une uniformisante de $\lambda$ dans $O_{\lambda}$ ($\pi_{\lambda}=\ell$ dans le cas non-ramifi\'e). Pour tout entier $n\geq 0$, les applications
\[ i_n : T_{\lambda}[\lambda^n]\rightarrow T_{\lambda}[\lambda^{n+1}], \ \ \ x\text{ mod }\lambda^nT_{\lambda}(A)\mapsto \pi_{\lambda}x\text{ mod }\lambda^{n+1}T_{\lambda}(A)\]
\noindent sont des morphismes, bien d\'efinis, injectifs. En prenant le syst\`eme inductif qu'ils forment, on note $T_{\lambda}[\lambda^{\infty}]$ la limite.

\medskip

\noindent Supposons dans la fin de ce paragraphe que $\ell$ est non ramifi\'e. Dans ce cas le $\O_{\ell}/\ell^n\O_{\ell}$-module $A[\ell^n]$ se d\'ecompose en le produit des $\O_{\lambda}/\lambda^n$-modules $T_{\lambda}[\lambda^n]$ dans le cas de type I, et, pour $\ell\notin S_{\text{ex}}(A)$, en le produit des $T_{\lambda}[\lambda^n]\oplus T_{\lambda}[\lambda^n]$ dans le cas de type II. De plus, par projection modulo $\ell^n$ (ou ce qui revient au m\^eme ici, modulo $\lambda^n$), on obtient :
\[\phi_{\lambda^{n}} : T_{\lambda}[\lambda^n]\times T_{\lambda}[\lambda^n]\rightarrow \O_{\lambda}/\ell^n\O_{\lambda}(1)\]
\noindent qui v\'erifie 
\[\phi_{\lambda^{n}}(\ell x,\ell y)=\phi_{\lambda^{n+1}}(x,y)^{\ell}.\]

\noindent Tout comme l'accouplement $\lambda$-adique, on voit par r\'eduction modulo $\ell^n$ que les accouplements $\phi_{\lambda^{n}}$ sont Galois \'equivariants, l'action de Galois \`a gauche se faisant via la repr\'esentation $\lambda$-adique, et \`a droite, via le caract\`ere cyclotomique $\ell$-adique usuel. 

\section{Modules isotropes}

\noindent On consid\`ere dans ce paragraphe une vari\'et\'e ab\'elienne $A/K$ g\'eom\'etriquement simple de type I ou II, telle que $\End_K(A)=\End_{\bar K}(A)$ et telle que $\End_K(A)$ est un ordre maximal de $\End_K(A)\otimes_{\Z}\Q$. Soit par ailleurs $\ell$ un premier et $\lambda$ une place au dessus de $\ell$. Notons $\pi_{\lambda}$ une uniformisante de $\lambda$. On suppose ici que $\ell$ est tel que la condition suivante (qui exclut un nombre fini de premiers) est r\'ealis\'ee : on a un accouplement bilin\'eaire altern\'e, non-d\'en\'en\'er\'e sur $T_{\lambda}(A)$ (qui est un $\O_{\lambda}$-module libre de rang $2h$) et sur le $\F_{\lambda}$-espace vectoriel $T_{\lambda}[\lambda]=T_{\lambda}(A)/\lambda T_{\lambda}(A)$.

\medskip

\defi Soit $H\subset T_{\lambda}[\lambda^{\infty}]$ un sous-groupe fini. Nous dirons que \textit{$H$ est totalement isotrope} si pour tous points $P,Q$ de $H\subset T_{\lambda}[\lambda^n]$, on a
\[\phi_{\lambda^n}(P,Q)=1,\]
\noindent o\`u $\phi_{\lambda^n}$ d\'esigne l'accouplement sur $T_{\lambda}[\lambda^n]$.

\medskip

\noindent Notons que si $H$ est totalement isotrope au sens pr\'ec\'edent, alors son sous-groupe des points de $\lambda$-torsion est totalement isotrope dans le $\F_{\lambda}$-espace vectoriel $T_{\lambda}[\lambda]$. On retrouve avec cette d\'efinition les deux lemmes suivants dont les preuves se reprennent mot pour mot du paragraphe 3.1 de \cite{HR} en rempla\c{c}ant $\Z$ par $\O_{\lambda}$ et $\ell$ par $\lambda$.

\begin{lemme}\label{h2}Soit $(e_1,\ldots,e_h)$ une base d'un sous-$\O_{\lambda}$-module isotrope maximal $H_{\infty}$ de $T_{\lambda}(A)$. Il existe un suppl\'ementaire $H_{\infty}'$ isotrope maximal et une base $(e_{h+1},\ldots,e_{2h})$ de celui-ci de sorte que dans la d\'ecomposition $T_{\lambda}(A)=H_{\infty}\oplus H'_{\infty}$ selon la base $(e_1,\ldots,e_{2h})$, la forme symplectique s'\'ecrit comme la forme canonique $J=\begin{pmatrix} 0 & I_h\cr -I_h&0\cr\end{pmatrix}$.
\end{lemme}

\begin{lemme}\label{hh2} Soit $n\geq 1$ et $H\subset T_{\lambda}[\lambda^{n}]$ un sous-groupe fini, totalement isotrope. Notons $\text{pr}_n : T_{\lambda}(A)\rightarrow T_{\lambda}[\lambda^n]$ la projection canonique modulo $\lambda^n$. Il existe un sous-groupe totalement isotrope $H_{\text{ti}}$ de $T_{\lambda}[\lambda^n]$, contenant $H$ et de m\^eme exposant et il existe un sous-$\O_{\lambda}$-module, $H_{\infty}$ de $T_{\lambda}(A)$, totalement isotrope, tel que $\text{pr}_n(H_{\infty})=H_{ti}$.
\end{lemme}

\rem Notons que dans \cite{HR} la version correspondante du lemme pr\'ec\'edent ne mentionne pas que l'on peut choisir $H_{ti}$ de m\^eme exposant que $H$. Toutefois la construction m\^eme de ce $H_{ti}$ fournie dans la preuve du lemme 3.7 de \cite{HR} donne imm\'ediatement cette information suppl\'ementaire.

\section{Propri\'et\'e $\mu$, version $\lambda$-adique}
\noindent On consid\`ere dans ce paragraphe une vari\'et\'e ab\'elienne $A/K$ g\'eom\'etriquement simple de type I ou II, telle que $\End_K(A)=\End_{\bar K}(A)$ et telle que $\End_K(A)$ est un ordre maximal de $\End_K(A)\otimes_{\Z}\Q$. On suppose par ailleurs ici que $\ell\notin S_{\text{ex}}(A)$. 

\subsection{Propri\'et\'e $\mu$}

\noindent \'Etant donn\'e un sous-groupe $H$ fini de $T_{\lambda}[\lambda^{\infty}]$, nous introduisons \`a pr\'esent l'invariant suivant~:
\[m_1(H)=\max\left\{k\in\N\ |\ \exists n\geq 0, \ \exists P,Q \in H\text{ d'ordre }\ell^n,\ \ \phi_{\lambda^n}(P,Q) \text{ est d'ordre }\ell^k\right\}.\]
\noindent Dire que $H$ est totalement isotrope \'equivaut \`a dire que $m_1(H)=0$. De plus on peut noter que, sur la d\'efinition, il est \'evident que $m_1(H)$ est sup\'erieur \`a la valeur $m$ suivante :
\[m(H):=\max\left\{k\in\N\ |\ \exists P,Q \in H\text{ d'ordre }\ell^k,\ \ \phi_{\lambda^k}(P,Q) \text{ est d'ordre }\ell^k\right\}.\]
\noindent Lorsque $H$ est de la forme $T_{\lambda}[\lambda^n]$, nous allons montrer que $m_1(T_{\lambda}[\lambda^n])=m(T_{\lambda}[\lambda^n])=n$.

\medskip

\defi Nous appelons \textit{propri\'et\'e $\mu$} pour une vari\'et\'e ab\'elienne le fait d'avoir, pour tout sous-groupe fini $H\subset T_{\lambda}[\lambda^{\infty}]$, l'\'egalit\'e \`a indice fini pr\`es, uniform\'ement en $(\ell,H)$~:
\[ K(\mu_{\ell^{m_1(H)}})\que K(H)\cap K(\mu_{\ell^{\infty}}).\]

\subsection{Propri\'et\'e $\mu$ pour $T_{\lambda}[\lambda^n]$}

\noindent Soit $n\geq 1$ un entier. La propri\'et\'e $\mu$ pour $T_{\lambda}[\lambda^n]$ d\'ecoule essentiellement formellement de la propri\'et\'e $\mu$ pour $A[\ell^n]$ et du fait que le multiplicateur $\mult_{\lambda}(\rho_{\lambda}(\sigma))$ est $\chi_{\ell}(\sigma)$. Plus pr\'ecis\'ement, on sait que concernant l'image de la representation $\lambda$-adique r\'esiduelle $\rho_{\lambda}$ (\`a valeur dans $\F_{\lambda}$), on a 


\begin{prop}\label{newprop} Soit $A$ de type I ou II et pleinement de type Lefschetz. On a les (presque) \'egalit\'es suivantes :

\begin{enumerate}
  \item $[\rho_{\ell}(G_K) , \rho_{\ell}(G_K)]\pre\prod_{\lambda | \ell} \Sp_{2h}(\F_{\lambda})\pre\Sp_{2h}(\O_E/\ell\O_E). $
  \item  $\rho_{\lambda}(G_K)\pre \{x\in\GSp_{2h}(\F_{\lambda})\ |\ \mult(x)\in\F_{\ell}^{\times}\}. $
\item $[\rho_{\ell^{\infty}}(G_K), \rho_{\ell^{\infty}}(G_K)] \pre \Hdg(A)(\Z_{\ell})=\prod_{\lambda\,|\,\ell} \Sp_{2h}(\O_{\lambda}). $
\item $\rho_{\lambda^{\infty}}(G_K)\pre\{x\in\GSp_{2h}(\O_{\lambda})\ |\ \mult(x)\in\Z_{\ell}^{\times}\}. $
\item 
$\rho_{\ell^{\infty}}(G_K)\pre \MT(A)(\Z_{\ell})=\{(x_{\lambda})\in\prod_{\lambda\,|\,\ell}\GSp_{2h}(\O_{\lambda})\ |\ \forall \lambda,\;\exists y \in\Z_{\ell}^{\times},\;\mult(x_{\lambda})=y\},$
\end{enumerate}

\noindent le produit portant sur les places $\lambda$ au dessus de $\ell$ dans l'anneau des entiers $\O_E$ de $E$.
\end{prop}

\demo  L'hypoth\`ese que $A$ est de type Lefschetz signifie que $\Hdg(A)=\Res_{E/\Q}\Sp_{E,2h}$ et $\MT(A)=\G_m\Res_{E/\Q}\Sp_{E,2h}$, l'hypoth\`ese que $A$ est pleinement de Lefschetz signifie que l'image de Galois est d'indice fini dans $\MT(A)(\Z_{\ell})$. Comme nous l'expliquons en appendice (th\'eor\`eme \ref{indicefini}), ceci entra\^{\i}ne en fait que cet indice est born\'e {\it ind\'ependamment} de $\ell$. En particulier l'indice de $\rho_{\lambda}(G_K)\cap\Sp_{2h}(\F_{\lambda})$ est born\'e    ind\'ependamment de $\ell$, disons par $c$. 
Observons maintenant  que $\Sp_{2h}(\F_{\ell^m})$ ne poss\`ede pas de sous-groupe d'indice ``petit" (ceci se voit en appliquant les lemmes 2.5 et 2.13 de \cite{HR}), c'est-\`a-dire que, pour $\ell\geq \ell_0=\ell_0({c})$, un sous-groupe d'indice inf\'erieur \`a $c$ est \'egal  au groupe  $\Sp_{2h}(\F_{\lambda})$ tout entier. 
D'apr\`es le lemme \ref{releveserre}  nous pouvons conclure que, pour $\ell\geq\ell_0$, nous avons $[\rho_{\lambda^{\infty}}(G_K),\rho_{\lambda^{\infty}}(G_K)]=\Sp_{2h}(\O_{\lambda})$. Ensuite en utilisant le fait que $\mult_{\lambda}(\rho_{\lambda^{\infty}}(\sigma))=\chi_{\ell^{\infty}}(\sigma)$ et que le caract\`ere cyclotomique est surjectif sur $\Z_{\ell}^{\times}$, toujours pour $\ell$ assez grand, on conclut que   
$\rho_{\ell^{\infty}}(G_K)\pre \MT(A)(\Z_{\ell})$, comme annonc\'e. Les autres \'egalit\'es s'en d\'eduisent ais\'ement. $\findemo$   

\medskip

\noindent On peut d\'eduire de ces consid\'erations l'observation suivante concernant la partie cyclotomique des extensions engendr\'ees, valable pour $A$ pleinement de type Lefschetz, de type I ou II. On peut d\'ecrire la situation via la tour d'extensions suivante :
\[
\xymatrix{
	K(A[\ell^{\infty}])\ar@{-}[dd]\ar@{-}[dr]	&																								\\
																						&	K(T_{\lambda}[\lambda^{\infty}]) \ar@{-}[dd]^{\Sp_{2h}(\O_{\lambda})}						\\
K(\mu_{\ell^{\infty}})\ar@{-}[dd]_{\Z_{\ell}^{\times}}\ar@{-}[dr]&																								\\				
																						& K(T_{\lambda}[\lambda^{\infty}])\cap K(\mu_{\ell^{\infty}})	\ar@{-}[dl]\\
												K										&																										}
 \]
\noindent Nous r\'esumons cela dans le corollaire suivant

\begin{cor}  Soit $A$ est de type I ou II, et pleinement de type Lefschetz.  On a les (presque) \'egalit\'es suivantes :
\[K(T_{\lambda}[\lambda^{\infty}])\cap K(\mu_{\ell^{\infty}})\pre K(\mu_{\ell^{\infty}}),\]
\noindent On a le m\^eme r\'esultat en niveau fini par r\'eduction modulo $\ell^n$.
\end{cor}
                                                
\subsection{Propri\'et\'e $\mu$ pour $H\subset T_{\lambda}[\lambda^n]$}

\begin{prop}\label{lmu} Soit $H$ un sous-groupe fini de $T_{\lambda}[\lambda^{\infty}]$. On a, uniform\'ement en $(\ell,H)$, l'in\'egalit\'e, 
\[[K(\mu_{\ell^{m_1(H)}}):\Q]\ll [K(H):\Q].\]
\end{prop}
\demo Soit $x,y\in H$ deux points de m\^eme ordre $\ell^n$ tels que $\phi_{\lambda^n}(x,y)$ est un \'el\'ement d'ordre $\ell^{m_1(H)}$. Montrons que l'extension $K(x,y)$ contient ``presque" $K(\mu_{\ell^{m_1(H)}})$. Ces deux extensions sont des sous-$K$-extensions de $K(T_{\lambda}[\lambda^n])$ et par la description du groupe de Galois de $K(T_{\lambda}[\lambda^n])/K$, on voit que le groupe de Galois $G_{x,y}$ de $K(T_{\lambda}[\lambda^n])$ sur $K(x,y)$ est donn\'ee par la presque \'egalit\'e suivante (valable pour tout $\ell$ assez grand), 
\[G_{x,y}\pre\{\rho_{\lambda^n}(\sigma)\in \GSp_{2h}(\O_{\lambda}/\ell^n\O_{\lambda})\ |\ \sigma\in G_K,\ \sigma \cdot x=x,\ \sigma\cdot y=y,\text{ et } \chi_{\ell^n}(\sigma)\in(\Z/\ell^n\Z)^{\times}\}.\]
\noindent Soit donc $\sigma\in G_K$ tel que $\rho_{\lambda^n}(\sigma)\in G_{x,y}$. On a
\[\phi_{\lambda^n}(x,y)=\phi_{\lambda^n}(\rho_{\lambda^n}(\sigma)(x),\rho_{\lambda^n}(\sigma)(y))=\chi_{\ell^n}(\sigma)\phi_{\lambda^n}(x,y).\]
\noindent On en d\'eduit que $\chi_{\ell^n}(\sigma)-1$ est un multiple de l'ordre de $\phi_{\lambda^n}(x,y)$ dans $\O_{\lambda}/\ell^n\O_{\lambda}$, autrement dit que $\chi_{\ell^n}(\sigma)=1\mod \ell^{m_1(H)}$. Or le groupe de Galois de $K(T_{\lambda}[\lambda^n])$ sur $K(\mu_{\ell^{m_1(H)}})$ est pr\'ecis\'ement constitu\'e des $\rho_{\lambda^n}(\sigma)$ tels que $\chi_{\ell^n}(\sigma)=1\mod \ell^{m_1(H)}$. On en d\'eduit le r\'esultat. $\findemo$

\medskip

\noindent Nous pouvons maintenant prouver la propri\'et\'e $\mu$ proprement dite~:

\begin{prop}\label{propmu} En notant $\delta(H):=\left(\Z_{\ell}^{\times}:\mult(G_0(H))\right)$ o\`u $G_0(H)=\Gal\left( K(A\left[{\lambda^{\infty}}\right])/K(H)\right)$, on a, pour tout sous-groupe $H$ fini de $T_{\lambda}[\lambda^{\infty}]$ l'\'egalit\'e \`a indice fini pr\`es, uniform\'ement en $(\ell,H)$,
\[ [K(H)\cap K(\mu_{\ell^{\infty}}):K]\que \delta(H).\]
\noindent De plus, pour tout $H$ sous-groupe fini de $T_{\lambda}[\lambda^{\infty}]$, on a l'inclusion suivante, qui est une \'egalit\'e \`a un indice fini pr\`es uniform\'ement en $(\ell,H)$~:
\[K(H)\cap K(\mu_{\ell^{\infty}}) \subset K(\mu_{\ell^{m_1(H)}}) \text{ et } K(H)\cap K(\mu_{\ell^{\infty}}) \que K(\mu_{\ell^{m_1(H)}}).\]
\end{prop}
\demo
\noindent On a la presque \'egalit\'e $\Gal\left( K(T_{\lambda}\left[{\lambda^{\infty}}\right])/K\right)\pre G_{\lambda}$ (introduit au lemme \ref{lemmeA}). Le groupe de Galois $\Gal\left( K(T_{\lambda}\left[{\lambda^{\infty}}\right])/K(\mu_{\ell^{\infty}})\right)$ s'identifie (c'est une presque \'egalit\'e) alors avec 
$SG_{\lambda}:=G_{\lambda}\cap \Ker(\mult)$.     Alors  $K(H)\cap K(\mu_{\ell^{\infty}})$ est la sous-extension fix\'ee par le groupe $U$ engendr\'e par $SG_{\lambda}$ et $G_0(H)$. On voit imm\'ediatement que le noyau de
$G_{\lambda}\buildrel{\mult}\over{\rightarrow}\Z_{\ell}^{\times}\rightarrow\Z_{\ell}^{\times}/\mult(G_0(H))$ est le groupe $U$ d'o\`u le premier \'enonc\'e.
\medskip

\noindent Passons maintenant \`a la seconde partie de la proposition. Commen\c{c}ons par consid\'erer $H_{\infty}$ un sous-groupe isotrope maximal de $T_{\lambda}(A)$. Par le lemme \ref{h2}, on peut supposer que dans une d\'ecomposition $T_{\lambda}(A)=H_{\infty}\oplus H'_{\infty}$ la forme symplectique s'\'ecrit comme la forme canonique $J$. On voit alors ais\'ement que,
\begin{align*}
\Gal\left(K(T_{\lambda}[\lambda^{\infty}])/K(H_{\infty})\right)        & \que\left\{M=\begin{pmatrix} I& *\cr 0& *\cr\end{pmatrix}\in\GSp_{2h}(\O_{\lambda})\ |\ \mult(M)\in \Z_{\ell}^{\times}\right\}\\
                                                                                                                                                                                                                & =\left\{M=\begin{pmatrix} I& S\cr 0& \alpha I\cr\end{pmatrix}\;|\; \alpha\in\Z_{\ell}^{\times}\;\text{et}\; S\;\text{sym\'etrique}\right\}.
\end{align*}
\noindent D'apr\`es le lemme \ref{lemmeA}, le groupe engendr\'e par ce dernier groupe et par le groupe
$\Sp_{2h}(\O_{\lambda})\pre\Gal\left(K(T_{\lambda}[\lambda^{\infty}])/K(\mu_{\ell^{\infty}})\right)$ est $\{x\in \GSp_{2h}(\O_{\lambda})\ |\ \mult(x)\in \Z_{\ell}^{\times}\}$ tout entier. Ainsi $K(H_{\infty})\cap K(\mu_{\ell^{\infty}})\que K$. Si $H$ est un sous-groupe fini de $T_{\lambda}[\lambda^{\infty}]$ totalement isotrope, dans ce cas le lemme \ref{hh2} et ce qui pr\'ec\`ede nous permettent de conclure :  on a $K(H)\cap K(\mu_{\ell^{\infty}})\que K$.

\medskip

\noindent Soit maintenant $H$ un sous-groupe fini non-isotrope de $T_{\lambda}[\lambda^{\infty}]$. Le groupe $[\ell^{m_1(H)}](H)$ est totalement isotrope. En effet si $P$ et $Q$ sont deux points d'ordre $\ell^n$ dans $H$, alors 
\[\phi_{\lambda^{n-m_1(H)}}(\ell^{m_1(H)}P,\ell^{m_1(H)}Q)=\phi_{\lambda^{n}}(P,Q)^{\ell^{m_1(H)}}=1 \text{ par d\'efinition de }m_1(H).\]
\noindent En appliquant le lemme \ref{hh2}, on trouve donc un sous-groupe $H'$ contenant $[\ell^{m_1}](H)$ de m\^eme exposant et il existe un sous-$\O_{\lambda}$-module $H_{\infty}$ totalement isotrope de $T_{\lambda}(A)$ tel que,  si, pour tout entier $n\geq 1$, $\text{pr}_n : T_{\lambda}(A)\rightarrow T_{\lambda}(A)/\ell^nT_{\lambda}(A)=T_{\lambda}[\lambda^n]$ d\'esigne la projection canonique, on a
\[\text{pr}_{r_H}(H_{\infty})=H'.\]
\noindent Par le lemme \ref{h2}, on peut supposer que dans une d\'ecomposition $T_{\lambda}(A)=H_{\infty}\oplus H'_{\infty}$ la forme symplectique s'\'ecrit comme la forme canonique $J$. Pour tout $n\geq 1$, notons 
\[H_n:=\text{pr}_n(H_{\infty})=H_{\infty}/H_{\infty}\cap\ell^n T_{\lambda}(A).\]
\noindent On a pour tout $n\geq 1$, $[\ell] H_{n+1}=H_n$. On peut donc poser
\[H^{\infty}=\bigcup_{n\geq 1}H_n\subset T_{\lambda}[\lambda^{\infty}].\]
\noindent De plus, on voit que, dans $K(T_{\lambda}[\lambda^{\infty}])$, le groupe de Galois correspondant \`a $H_{\infty}$ est le m\^eme que celui correspondant \`a $H^{\infty}$. On a 
\[H\subset[\ell^{m_1(H)}]^{-1}(H')=[\ell^{m_1(H)}]^{-1}(H_{r_H})\subset[\ell^{m_1(H)}]^{-1}(H^{\infty}).\]
\noindent En consid\'erant la multiplication par $[\ell^{m_1(H)}]$ sur $H^{\infty}$, on en d\'eduit (car $H^{\infty}$ est $\ell$-divisible)   que 
\[H\subset H^{\infty}+\ker[\ell^{m_1(H)}]=:\widetilde{H^{\infty}},\]
\noindent o\`u $[\ell^n]$ est le morphisme de multiplication dans $T_{\lambda}[\lambda^{\infty}]$. Ainsi comme dans le cas totalement isotrope, on se ram\`ene \`a une situation o\`u un lemme de groupe permet de conclure : le groupe de Galois $\Gal\left(K(T_{\lambda}[\lambda^{\infty}])/K(\widetilde{H^{\infty}})\right)$ n'est autre que (il s'agit d'une \'egalit\'e $\que$ \`a indice fini pr\`es) 
\[\left\{M\in\GSp_{2h}(\O_{\lambda})\;|\;\forall i\leq g\ Me_{g+i}=e_{g+i}\text{ mod }\ell^{m_1(H)},\ Me_i=e_i\text{ et }\mult(M)\in\Z_{\ell}^{\times}\right\}.\]
\noindent La m\^eme preuve  que celle du corollaire 2.11 de \cite{HR} donne alors le r\'esultat : le groupe engendr\'e par $\Gal\left(K(T_{\lambda}[\lambda^{\infty}])/K(\widetilde{H^{\infty}})\right)$ et $\Sp_{2h}(\O_{\lambda})$ est (avec une \'egalit\'e $\que$ \`a indice fini pr\`es) 
\[\left\{M\in\GSp_{2h}(\O_{\lambda})\ |\ \mult(M)\in \Z_{\ell}^{\times}\text{ et } \mult(M)\equiv 1\mod\ell^{m_1(H)}\right\}.\]
\noindent Notamment, on a,
\[K(H)\cap K(\mu_{\ell^{\infty}})\subset K(\widetilde{H^{\infty}})\cap K(\mu_{\ell^{\infty}})\subset K(\mu_{\ell^{m_1(H)}}),\]
\noindent la seconde inclusion \'etant aussi une \'egalit\'e \`a indice fini pr\`es, ie $K(\widetilde{H^{\infty}})\cap K(\mu_{\ell^{\infty}})\que K(\mu_{\ell^{m_1(H)}})$. La proposition pr\'ec\'edente \ref{lmu} permet de conclure. $\findemo$

\section{Preuve du th\'eor\`eme principal pour $H\subset A[\ell]$\label{p6}}

\noindent Soit $A/K$ une vari\'et\'e ab\'elienne sur un corps de nombres, telle que $\End_K(A)=\End_{\bar K}(A)$. On commence par se ramener au cas $\ell$-adique (cf. \cite{hindry-ratazzi1} proposition 4.1) grace \`a la presque ind\'ependance rappel\'ee \`a la proposition \ref{indep}~: 
\begin{prop}\label{ladique}Soit $\alpha>0$. Pour d\'emontrer que $\gamma(A)\leq\alpha$, il suffit de montrer que : il existe une cons\-tan\-te strictement positive $C(A/K)$ ne d\'ependant que de $A/K$ telle que pour tout nombre premier $\ell$, pour tout sous-groupe fini $H_{\ell}$ de $A[\ell^{\infty}]$, on a
\begin{equation}\label{equa}
\text{Card}\left(H_{\ell}\right)\leq C(A/K)[K(H_{\ell}):K]^{\alpha}.
\end{equation}
\end{prop}

\rem Rappelons que l'on a suppos\'e que la vari\'et\'e ab\'elienne $A/K$ est telle que $\End_K(A)=\End_{\bar K}(A)$. Concernant notre question de borne sur la torsion, ceci nous permet de supposer que le groupe fini $H\subset A[\ell^n]$ est en fait un $\End_{\bar K}(A)$-module. En effet : notons $H_E$ le $\End_{\bar K}(A)$-module engendr\'e par $H$ et supposons que l'on ait pour $H_E$ une in\'egalit\'e de la forme suivante, uniform\'ement en $(\ell,H)$,
\[|H_E|\ll [K(H_E):K]^{\alpha}.\]
\noindent on a donc $|H|\ll [K(H_E):K]^{\alpha}$ car $H$ est inclus dans $H_E$. Mais $\End_K(A)=\End_{\bar K}(A)$, donc si $x\in H$ et $f\in\End_{\bar K}(A)$ alors $f(x)$ est encore un point de $A$ qui est $K(H)$ rationnel, donc $K(H)=K(H_E)$. En particulier ceci implique que $|H|\ll [K(H):K]^{\alpha}$ comme annonc\'e.

\medskip

\noindent Nous nous pla\c{c}ons dans toute la suite de ce paragraphe dans la situation particuli\`ere d'une vari\'et\'e ab\'elienne $A$  d\'efinie sur $K$, g\'eom\'etriquement simple de type I ou II, qui est pleinement de type Lefschetz. Nous supposerons de plus que $\ell\notin S_{\text{ex}}(A)$ de sorte \`a pouvoir appliquer les techniques d\'evelopp\'ees dans les paragraphes pr\'ec\'edents. Enfin nous prenons le cas particulier d'une situation horizontale d'un sous-groupe $H$ de $A[\ell]$ (en particulier il s'agit d'un $\F_{\ell}$-espace vectoriel). Par la remarque pr\'ec\'edente, nous pouvons m\^eme supposer que $H$ est un $\O_E/\ell\O_E$-module. Nous avons la d\'ecompo\-sition suivante : 
\[H=\left\{\begin{matrix}\prod_{\lambda|\ell} H[\lambda]\subset \prod_{\lambda|\ell}T_{\lambda}[\lambda]& (\text{Type I})\cr \prod_{\lambda|\ell} H[\lambda]\oplus H[\lambda]\subset \prod_{\lambda|\ell}T_{\lambda}[\lambda]\oplus T_{\lambda}[\lambda] & (\text{Type II})\cr\end{matrix}\right..\]

\noindent On sait par la proposition \ref{newprop} que pour tout $\ell$ on a,
\[\rho_{\lambda}(G_K)\pre\left\{M\in\GSp_{2h}(\F_{\lambda})\ |\ \mult(M)\in \F_{\ell}^{\times}\right\}.\] 

\noindent Dans notre situation les $H[\lambda]\subset T_{\lambda}[\lambda]$ sont des $\F_{\lambda}$-espaces vectoriels. Rappelons que l'on a, uniform\'ement en $(\ell,H)$, l'\'egalit\'e \`a indice fini pr\`es $\delta(H[\lambda])\que [K(H[\lambda])\cap K(\mu_{\ell}):K]$. On obtient ainsi : 

\begin{lemme}\label{muH}Si $H[\lambda]$ est inclus dans un sous-espace totalement isotrope du $\F_{\lambda}$-ev $T_{\lambda}[\lambda]$ alors, uniform\'ement en $(\ell,H)$, on a $\delta(H[\lambda])\que 1$. Sinon $\delta(H[\lambda])\que \ell$.
\end{lemme}

\medskip

\begin{lemme}\label{pmu1} \label{calculH1} Uniform\'ement en $(\ell,H)$, on a 
\[\delta(H[\lambda])\que\left(\F_{\ell}^{\times}:\mult(G_0(H[\lambda]))\right)\text{ o\`u }\ G_0(H[\lambda])=\Gal\left( K(T_{\lambda}[\lambda])/K(H[\lambda])\right).\]
\noindent On a de plus~:
\[[K(H[\lambda]):K]\pre(\rho_{\lambda}(G_K):G_0(H[\lambda]))\pre\delta(H[\lambda])(\Sp_{2h}[\F_{\lambda}):G(H[\lambda])).\]
\end{lemme}
\demo Pour le premier point, on a $\Gal\left(K(T_{\lambda}[\lambda])/K\right)\pre\rho_{\lambda}(G_K)$. Le groupe de Galois $\Gal\left(K(T_{\lambda}[\lambda])/K(\mu_{\ell})\right)$ est alors presque \'egal \`a $SG_{\lambda}:=\rho_{\lambda}(G_K)\cap \Ker(\mult)$. Alors  $K(H[\lambda])\cap K(\mu_{\ell})$ est la sous-extension fix\'ee par le groupe $U$ engendr\'e par $SG_{\lambda}$ et $G_0(H[\lambda])$. On voit imm\'ediatement que le noyau de
$\rho_{\lambda}(G_K)\buildrel{\mult}\over{\rightarrow}\F_{\ell}^{\times}\rightarrow\F_{\ell}^{\times}/\mult(G_0(H[\lambda]))$ est le groupe $U$. Pour le second point : la premi\`ere \'egalit\'e est donn\'ee par la th\'eorie de Galois car on a que $\Gal(K(T_{\lambda}[\lambda])/K)\pre\rho_{\lambda}(G_K)$. La seconde \'egalit\'e est une chasse au diagramme facile.\hfill$\Box$

\medskip

\noindent Notons $d_{\lambda}$ la dimension de $H[\lambda]$ sur $\F_{\lambda}$ et $(e_1,\ldots e_{d_{\lambda}})$ une base que l'on compl\`ete en une base $(e_1,\ldots,e_{2h})$ de $T_{\lambda}[\lambda]$. On d\'efinit
\[G(H[\lambda])=\left\{M\in\Sp_{2h}(\F_{\lambda})\;|\; Me_i= e_i,\ 1\leq i\leq d_{\lambda}\right\}.\]

\noindent Notons $(\hat{e}_1,\ldots,\hat{e}_{2h})$ une base de $T_{\lambda}(A)$ relevant la base sur $\F_{\lambda}$. Introduisons maintenant le groupe alg\'ebrique sur $\O_{\lambda}$ suivant :
\[ G_{1}:=\left\{M\in \Sp_{2h}\ |\ M\hat{e}_i=\hat{e}_i,\ 1\leq i\leq d_{\lambda}\right\}.\]
\noindent On voit que
\[G(H[\lambda])=\left\{M\in \Sp_{2h}(\F_{\lambda})\ |\  M\in G_1\mod\lambda\right\}.\]
\noindent Par changement de base symplectique sur $\F_{\lambda}$, $G_1$ est conjugu\'e sur $\F_{\lambda}$ \`a l'un des groupes $P_{r,s}$ introduits pr\'ec\'edemment. En posant $G=\Sp_{2h}$, et en rappelant que $\text{Card}(\F_{\lambda})=\ell^{f(\lambda)}$, on voit que, d'apr\`es lemme \ref{cle} on a
\[ [K(H[\lambda]):K]\que \ell^{m(H[\lambda])}\ell^{f(\lambda)\codim P_{r_{\lambda},s_{\lambda}}},\]
\noindent o\`u $(r_{\lambda},s_{\lambda})$ (avec \'eventuellement $s_{\lambda}=0$) est le couple correspondant \`a $H[\lambda]$. 

\medskip

\noindent Utilisant le lemme \ref{muH}, le th\'eor\`eme 6.6 de \cite{hindry-ratazzi1} s'adapte imm\'ediatement (cf. la proposition \ref{refauprodce} ci-apr\`es) pour donner :

\begin{prop} Avec les notations pr\'ec\'edentes, uniform\'ement en $(\ell,H)$, on a 
\[ \ell^{m(H)}\que[K(H)\cap K(\mu_{\ell}):K]\que\max_{\lambda |\ell} \ell^{m(H[\lambda])}\]
\noindent et
\[ [K(H):K(\mu_{\ell^{m(H)}})]\que\prod_{\lambda |\ell}\left[K(H[\lambda]):K(\mu_{\ell^{m(H[\lambda])}})\right].\]
\end{prop}

\subsection{Cas totalement d\'ecompos\'e}
\noindent Nous supposons ici que $\ell\notin S_{\text{ex}}(A)$ est totalement d\'ecompos\'e dans $\O_E$. Notre situation est alors la suivante : 
\[H=\left\{\begin{matrix}\prod_{\lambda|\ell} H[\lambda]\subset \prod_{\lambda|\ell}T_{\lambda}[\lambda] & \text{ et } & \rho_{\ell}=\prod_{\lambda |\ell}\rho_{\lambda}& (\text{Type I})\cr \prod_{\lambda|\ell} H[\lambda]\oplus H[\lambda]\subset \prod_{\lambda|\ell}T_{\lambda}[\lambda]\oplus T_{\lambda}[\lambda] & \text{ et } & \rho_{\ell}=\prod_{\lambda |\ell}\rho_{\lambda}\oplus\rho_{\lambda}& (\text{Type II})\cr\end{matrix}\right..\]
\noindent De plus 
\[\Gal(K(A[\ell])/K(\mu_{\ell}))=\prod_{\lambda |\ell} \Gal(K(T_{\lambda}[\lambda])/K(\mu_{\ell})).\]
\noindent Du point du vue combinatoire, les formules sont identiques \`a celles d'un     produit de vari\'et\'es ab\'eliennes de type $\GSp_{2h}$ et, les r\'esultats du paragraphe pr\'ec\'edent nous indiquent que la combinatoire n'est finalement autre que celle d'un produit de $e$ vari\'et\'es ab\'eliennes de type $\GSp_{2h}$, deux \`a deux non-isog\`enes. Nous pouvons donc directement en d\'eduire la valeur de l'exposant $\gamma(A)$.

\medskip

\defi Nous noterons dans la suite : $d=1$ si $A$ et de type I et $d=2$ si A est de type II.

\medskip

\noindent Les calculs de $\cite{HR}$ (paragraphes 4.1 et 6.2) donnent dans ce cas :
\[\gamma(A)=\sup_{I\subset\{1,\ldots,e\}}\frac{2\sum_{\lambda \in I}dh}{1+(2h^2+h)|I|}.\]
\noindent Ce sup se calcule ais\'ement (le max est atteint pour $I=\{1,\ldots,e\}$) et on trouve donc
\[\gamma(A)=\frac{2dhe}{1+(2h^2+h)e}=\frac{2\dim A}{1+\dim \Res_{E/\Q}\Sp_{2h}}=\frac{2\dim A}{\dim \MT(A)}.\]

\subsection{Cas g\'en\'eral}
\noindent Nous ne supposons plus d\'esormais que $\ell$ est totalement d\'ecompos\'e, la combinatoire qui r\'esulte est donc diff\'erente et il faut dans ce cadre g\'en\'eral la refaire explicitement (ceci contient d'ailleurs le cas du sous-paragraphe pr\'ec\'edent). On a

\[H[\lambda]=\left(\O_{\lambda}/\lambda\right)^{r_{\lambda}+s_{\lambda}}\ \text{ avec }\ s_{\lambda}=0 \text{ ssi $H[\lambda]$ est inclus dans un Lagrangien}.\]
\noindent De plus on a, quitte \`a r\'eordonner, 
\[0\leq s_{\lambda}\leq r_{\lambda}\leq h \text{ o\`u } 2h=\dim_{\F_{\lambda}}T_{\lambda}[\lambda], \text{ et } \sum_{\lambda | \ell }f(\lambda)=[E:\Q]=e\ \text{ et }\ dhe=g=\dim A,\]
\noindent o\`u l'on note comme pr\'ec\'edemment $d=1$ si $A$ est de type I et $d=2$ si $A$ est de type II. 

\medskip

\noindent  On obtient finalement, sous les conditions pr\'ec\'edentes, la valeur suivante pour le cardinal de $H$ :
\[\text{Card}(H)=\ell^{d\sum_{\lambda | \ell}f(\lambda)(r_{\lambda}+s_{\lambda})}.\]
\noindent Le degr\'e de l'extension $[K(H):K]$ d\'epend selon que les $H[\lambda]$ sont ou non inclus dans des Lagrangiens. Si l'un des $H[\lambda]$ n'est pas inclus dans un Lagrangien alors nous obtenons
\[[K(H):K]\que\ell^{1+\sum_{\lambda | \ell}f(\lambda)\codim P_{r_{\lambda},s_{\lambda}}}.\]
\noindent Si par contre tout les $H[\lambda]$ sont inclus dans un Lagrangien alors nous obtenons
\[[K(H):K]\que\ell^{\sum_{\lambda | \ell}f(\lambda)\codim P_{r_{\lambda},0}}.\]
\noindent Il reste \`a interpr\'eter le quotient de l'exposant de $\ell$ du $\text{Card(H)}$ par celui de $[K(H):K]$ pour conclure : c'est l'objet du paragraphe combinatoire suivant.

\subsection{Combinatoire}

\noindent Comme pr\'ec\'edemment on note $d=1$ si $A$ est de type I et $d=2$ si $A$ est de type II. Nous sommes ramen\'es \`a calculer la quantit\'e:

\[\frac{1}{d}\gamma:=\max_{r_{\lambda},s_{\lambda}}\frac{\sum_{\lambda\,|\,\ell}f(\lambda)(r_{\lambda}+s_{\lambda})}
{\delta+\sum_{\lambda\,|\,\ell}f(\lambda)\codim P_{r_{\lambda},s_{\lambda}}}\]
\noindent o\`u le maximum est pris pour  $0\leq s_{\lambda}\leq r_{\lambda}\leq h$ et $\delta$ vaut $0$ (resp. $1$) si tous les $s_{\lambda}$ sont nuls (resp. si l'un des $s_{\lambda}$ est non nul).

\begin{prop} Soit $\gamma=\gamma(A)$ d\'efini ci-dessus, alors
\[\gamma=\frac{2dhe}{1+2eh^2+he}=\frac{2\dim A}{\dim\MT(A)}.\]
\end{prop}

\medskip

\noindent Nous donnons ci-dessous, dans le cas particulier de la proposition ci-dessus, une preuve via les interpolateurs de Lagrange. Un argument combinatoire diff\'erent, sera donn\'e plus loin dans le cas g\'en\'eral de la preuve du lemmme \ref{comb}, l'argument suivant n'est donc pas indispensable mais a l'avantage d'\^etre assez direct.  

\noindent{\bf Remarques ``num\'eriques". } 
\begin{enumerate}
\item On peut r\'e\'ecrire, pour $P_{r,s}\subset \Sp_{2g}$:
$$\codim P_{r,s}=2g(r+s)-rs-\left(\frac{r^2+s^2}{2}\right)+\frac{r+s}{2}=\left(2g+\frac{1}{2}-\frac{r+s}{2}\right)(r+s).$$
\noindent On observe en particulier que la dimension ou codimension de $P_{r,s}$ ne d\'epend que de $r+s$.
\item Nous allons devoir \'etudier le sens de variation d'une fraction rationnelle de la forme:
$$f(x)=\frac{a-x}{A-\frac{x(x-1)}{2}}$$
\noindent dont la d\'eriv\'ee s'\'ecrit:
$$f'(x)=-2\frac{(x-a)^2+2A+a-a^2}{(2A-x(x-1))^2}$$
\noindent et est donc d\'ecroissante d\`es que $2A+a\geq a^2$.
\end{enumerate}

\demo La preuve consiste \`a appliquer le calcul diff\'erentiel \`a la fonction de variables $(\underline{r},\underline{s}):=(r_{\lambda},s_{\lambda})_{\lambda | \ell}$ dont on veut \'evaluer le maximum:
$$ \psi(\underline{r},\underline{s}):=\frac{N(\underline{r},\underline{s})}{D(\underline{r},\underline{s})}:=\frac{\sum_{\lambda\,|\,\ell}f(\lambda)(r_{\lambda}+s_{\lambda})}
{\delta+\sum_{\lambda\,|\,\ell}f(\lambda)\codim P_{r_{\lambda},s_{\lambda}}}
$$
(nous \'ecrivons la fonction sous la  forme 
``Num\'erateur/D\'enominateur= N/D").
Commen\c{c}ons par traiter le cas o\`u tous les $s_{\lambda}$ sont nuls. Les diff\'erentielles des deux fonctions $N$ et $D$ s'\'ecrivent
$$\partial N=(f(\lambda))_{\lambda |\ell}\qquad{\rm et}\qquad \partial D=\left(f(\lambda)(2h-r_{\lambda}+\frac{1}{2})  \right)_{\lambda |\ell}$$
Le th\'eor\`eme de Lagrange indique que, en un  maximum de $N/D$, ces deux diff\'erentielles sont proportionnelles, donc $2h-r_{\lambda}+\frac{1}{2}$ est constant, ou encore, $r_{\lambda}=2h-\kappa$ (avec $h\leq\kappa\leq 2h$). On obtient alors
$$\frac{N}{D}=\frac{2he-\kappa e}{2\sum_{\lambda|\ell} f(\lambda)h^2+he-e\kappa(\kappa-1)/2}=\frac{2h-\kappa}{2h^2+h-\kappa(\kappa-1)/2}$$
La fonction \`a droite est d\'ecroissante avec $\kappa$ et donc major\'ee par la valeur en $\kappa=h$, c'est-\`a-dire
$2/(3h+3)$ (noter que $\kappa\geq h$).

On traite ensuite le cas g\'en\'eral (avec l'un des $s_{\lambda}$ non nul), on pose donc
$$N=:\sum_{\lambda\,|\,\ell} f(\lambda)(r_{\lambda}+s_{\lambda});\quad D:=1+\sum_{\lambda\,|\,\ell} f(\lambda)(r_{\lambda}+s_{\lambda})\left(2h+\frac{1}{2}-\frac{r_{\lambda}+s_{\lambda}}{2}\right)$$
Le th\'eor\`eme de Lagrange indique maintenant que, en un maximum  de $N/D$, on aura $s_{\lambda}+r_{\lambda}=2h-\kappa$, avec maintenant $0\leq\kappa\leq 2h$. En reportant on obtient:
$$\frac{N}{D}\leq \frac{e(2h-\kappa)}{1+2\sum_{\lambda\,|\,\ell} f(\lambda)h^2+eh-e\kappa(\kappa-1)/2}= \frac{2h-\kappa}{\frac{1}{e}+2h^2+h-\kappa(\kappa-1)/2}$$
Cette derni\`ere fonction est d\'ecroissante en $\kappa$ donc major\'ee par la valeur en $\kappa=0$, ce qui donne au final:
$$\psi\leq \max\left\{ \frac{2}{3(h+1)}, \frac{2he}{1+2eh^2+eh}\right\}=\frac{2he}{1+2eh^2+eh}$$
Observons que $\kappa=0$ corre\'eond \`a $r_{\lambda}+s_{\lambda}=2h$ donc \`a $r_{\lambda}=s_{\lambda}=h$.
En consid\'erant donc  le cas $r_{\lambda}=s_{\lambda}=h$, on obtient
$$\psi=\frac{2he}{1+2h^2e+he}.$$ 

\section{Cas d'un groupe $H$ quelconque\label{p7}}
\noindent Dans ce paragraphe nous allons donner une preuve du r\'esultat principal (th\'eor\`eme \ref{theoproduit}). Rappelons que l'on a suppos\'e que la vari\'et\'e ab\'elienne $A/K$ est un produit $\prod_{i=1}^dA_i^{n_i}$ de vari\'et\'es ab\'eliennes simples, chacune de type I ou II et chacune pleinement de type Lefschetz. Nous avons d\'ej\`a indiqu\'e que l'on peut supposer de plus que pour tout $i$, les $A_i$ sont telles que $\End_{\bar K}(A_i)=\End_K(A_i)$ et telles que $\End_K(A_i)$ est un ordre maximal dans $\End_K(A_i)\otimes_{\Z}\Q$.

\medskip

\defi \label{ex2} Avec la notation de la d\'efinition 3.4, nous noterons dans la suite $S_{\text{ex}}=\bigcup_{i=1}^dS_{\text{ex}}(A_i)$.

\medskip 
 
\noindent Dans la suite de ce paragraphe nous supposerons que $\ell\notin S_{\text{ex}}$.

\medskip

\noindent Soit $H$ un sous-groupe fini de $A[\ell^{\infty}]$. Par le paragraphe 4.2 de \cite{hindry-ratazzi1}, on peut supposer que $H$ s'\'ecrit sous la forme $H=\prod_{i=1}^dH_i^{n_i}$. De plus par la remarque du paragraphe pr\'ec\'edent nous pouvons supposer que chaque $H_i$ est un $\O_{\ell,i}$-module, inclus dans un $A_i[\ell^n]$ pour $n$ convenable (o\`u l'on note $\O_{\ell,i}$ le tensoris\'e par $\Z_{\ell}$ de $\End(A_i)$). Notons

\[\I_{\ell}:=\left\{(\lambda,i)\ |\ i\in\{1,\ldots,d\},\ \text{ et }\lambda \text{ une place du centre de $\End(A_i)\otimes\Q$ au-dessus de $\ell$}\right\}.\]
\noindent Pour $(\lambda,i)\in \I_{\ell}$, posons $\O_{\lambda,i}$ la composante $\lambda$-adique de $\O_{\ell,i}$ et posons $X_{\lambda,i}$ le morceau de $H_i$ correspondant \`a $\lambda$. Dans le cas de type II, $X_{\lambda,i}$ se d\'ecompose \`a son tour en deux copies isomorphes : $X_{\lambda,i}=H_{\lambda,i}\oplus H_{\lambda,i}$. Dans le cas de type I, on pose $H_{\lambda,i}:=X_{\lambda,i}$. Avec des notations \'evidentes, les $H_{\lambda,i}$ sont des $\O_{\lambda,i}/\ell^n\O_{\lambda,i}$-sous-modules de $T_{\lambda,i}[\lambda^n]$. On a finalement la d\'ecomposition suivante de $H$ :
\[H=\prod_{(\lambda,i)\in \I_{\ell}}X^{n_i}_{\lambda,i}.\]
\noindent En tant que groupe on sait que pour $(\lambda,i)\in \I_{\ell}$, 
\[\O_{\lambda,i}/\ell^n\O_{\lambda,i}=(\Z/\ell^n\Z)^{f(\lambda)} \text{ et } T_{\lambda,i}[\lambda^n]=(\Z/\ell^n\Z)^{2h_if(\lambda)}.\]
\noindent  On sait \'egalement que, uniform\'ement en   $\ell$ et  $\lambda$, on a,
\[\rho_{A_i,\lambda^{\infty}}(G_K)\pre\left\{M\in\GSp_{2h_i}(\O_{\lambda,i})\ |\ \mult(M)\in \Z_{\ell}^{\times}\right\}.\] 
\noindent Dans tous les cas on obtient ainsi une \'egalit\'e \`a indice fini pr\`es, en r\'eduisant modulo $\ell^n$.

\medskip

\noindent Soit donc $H_{\lambda,i}$ un sous-groupe fini de $T_{\lambda,i}[\lambda^{\infty}]$, on pose
\[G_0(H_{\lambda,i}):=\left\{M\in\GSp_{2h_i}(\O_{\lambda,i})\ |\ \mult(M)\in\Z_{\ell}^{\times},\ \forall x\in H_{\lambda,i},\;Mx=x\right\}.\]
\noindent et $G(H_{\lambda,i}):=G_0(H_{\lambda,i})\cap\Sp_{2h_i}(\O_{\lambda,i})$. Comme $\O_{\lambda,i}/\ell^n\O_{\lambda,i}$-module et comme groupe abstrait, $H_{\lambda,i}$ est de la forme 
\[H_{\lambda,i}\simeq\prod_{j=1}^{2h_i}\O_{\lambda,i}/\ell^{m_j}\O_{\lambda,i}\simeq\prod_{j=1}^{2h_i}(\Z/\ell^{m_j}\Z)^{f(\lambda)},\]
\noindent o\`u nous sous-entendons, pour ne pas alourdir plus que de raison les notations, que les nombres $m_j$ d\'ependent \'egalement de $(\lambda,i)$.

\medskip

\noindent  Notons $e_{(\lambda,i)}^1, \ldots, e_{(\lambda,i)}^{2h_i}$ un syst\`eme de g\'en\'erateurs (en tant que $\O_{\lambda,i}/\ell^n\O_{\lambda,i}$-module) ; les $e_{(\lambda,i)}^j$ \'etant d'ordre respectifs $\ell^{m_{j}}$. Notons de plus $\{\hat{e}_{(\lambda,i)}^1,\ldots,\hat{e}_{(\lambda,i)}^{2h_i}\}$ une base du $\O_{\lambda,i}$-module libre $T_{\lambda,i}:=T_{\lambda}(A_i)$ relevant la famille $\{e_{(\lambda,i)}^j\}$, \textit{i.e.} telle que $e_{(\lambda,i)}^j=\hat{e}_{(\lambda,i)}^j\mod \ell^{m_j}$ pour tout $j$. On a
\[G(H_{\lambda,i})=\left\{M\in\Sp_{2h_i}(\O_{\lambda,i})\;|\; M\hat{e}_{(\lambda,i)}^j= \hat{e}_{(\lambda,i)}^j\mod\ell^{m_j},\ 1\leq j\leq 2h_i\right\}.\]

\medskip

\begin{lemme} \label{calculH1a} Notons $\delta(H_{\lambda,i}):=\left(\Z_{\ell}^{\times}:\mult(G_0(H_{\lambda,i}))\right)$. Uniform\'ement en $(\ell,H)$, on a alors~:
\[[K(H_{\lambda,i}):K]\que(\rho_{A_i,\lambda^{\infty}}(G_K):G_0(H_{\lambda,i}))\que\delta(H_{\lambda,i})(\Sp_{2h_i}(\O_{\lambda,i}):G(H_{\lambda,i})).\]
\end{lemme}
\demo Comme le lemme \ref{calculH1}. $\findemo$

\medskip

\noindent Quitte \`a renum\'eroter on peut supposer que les exposants $m_j$ (correspondants aux $e_{(\lambda,i)}^j$) sont ordonn\'es dans l'ordre d\'ecroissant : $m_1\geq\ldots\geq m_{2h_i}$. On pose alors
\[m^1:=\max \{m_i\ |\ m_i\not=0\}\ \text{ et par r\'ecurrence}\ m^{r+1}=\max\{m_i\ |\ m_i<m^r\}.\]
\noindent On obtient ainsi une suite strictement d\'ecroissante $m^1>\ldots > m^{t_{\lambda,i}}\geq 1$ (avec $t_{\lambda,i}\leq 2h_i$). Le groupe $H_{\lambda,i}$ est isomorphe \`a $\prod_{j=1}^{t_{\lambda,i}}\left(\Z/\ell^{m^j}\Z\right)^{f(\lambda)a_j}$, les $a_j$ d\'ependants de $(\lambda,i)$. On d\'efinit ensuite pour tout $1\leq r\leq t_{\lambda,i}$, les sous-ensembles emboit\'es
\[ I_{r}=\{j\in\{1,\ldots,2h_i\}\ |\ m_j\geq m^{r}\}\ \ \text{ de cardinal }\ \ \left|I_{r}\right|=\sum_{j=1}^{r}a_j.\]
\noindent Introduisons maintenant la suite croissante de groupes alg\'ebriques sur $\O_{\lambda,i}$ suivants :
\[\forall 1\leq r\leq t_{\lambda,i}\ \  G_{r,(\lambda,i)}:=\left\{M\in \Sp_{2h_i}\ |\ M\hat{e}_{(\lambda,i)}^j=\hat{e}_{(\lambda,i)}^j\ \ \forall j\in I_{t_{\lambda,i}+1-r}\right\}.\]
\noindent On voit que
\[G(H_{\lambda,i})=\left\{M\in \Sp_{2h_i}(\O_{\lambda,i})\ |\ \forall 1\leq r\leq t_{\lambda,i}\ \text{on a, } M\in G_{r,(\lambda,i)}\mod\ell^{m^{t_{\lambda,i}+1-r}}\right\}.\]
\noindent Par changement de base symplectique sur $\F_{\lambda}$, le couple $(\lambda,i)$ \'etant fix\'e, chacun des $G_{j,(\lambda,i)}$ est conjugu\'e sur $\F_{\lambda}$ \`a l'un des groupes $P_{r,s}$ introduits au paragraphe \ref{groupe}. En posant $G=\Sp_{2h_i}$, on voit que, avec les notations du lemme \ref{cle}, on a 
\[G(H_{\lambda,i})=H(m^1,\ldots,m^{t_{\lambda,i}}).\]
\noindent On va donc pouvoir appliquer le lemme \ref{cle}.

\subsection{Cas d'un morceau $H_{\lambda,i}$\label{m1}}

\noindent 
Le couple $(\lambda,i)$ \'etant fix\'e nous renoterons dans ce paragraphe $t:=t_{\lambda,i}$ afin de soulager un peu les notations. On peut appliquer le lemme \ref{cle}, uniform\'ement en $(\ell,H)$, on a~:
\[\left(\Sp_{2h_i}(\O_{\lambda,i}):G(H_{\lambda,i})\right)\gg \ell^{\sum_{j=1}^{t}f(\lambda)\codim(G_{j,(\lambda,i)})(m^{t+1-j}-m^{t+1-(j-1)})},\]
\noindent o\`u l'on a pos\'e $m^{t+1}=0$ et o\`u $\codim(G_{j,(\lambda,i)})$ est la codimension de $G_{j,(\lambda,i)}$ dans $\Sp_{2h_j}$. Les groupes alg\'ebriques $G_{j,(\lambda,i)}$ \'etant conjugu\'es sur $\F_{\lambda}$ aux $P_{r,s}$ (avec \'eventuellement $s=0$), $\codim(G_{j,(\lambda,i)})$ est \'egalement la codimension du groupe $P_{r_j,s_j}$ correspondant. Par ailleurs, la suite des $(G_{j,(\lambda,i)})_j$ \'etant croissante, $(\lambda,i)$ \'etant fix\'e, la suite des $(P_{r_j,s_j})_j$ l'est \'egalement. Ceci se traduit par
\[\forall j,\ \ \ r_j\geq r_{j+1}\ \ \ \text{ et }\ \ \ s_j\geq s_{j+1}.\]
\noindent Il nous reste \`a calculer la valeur de $\delta(H_{\lambda,i})$ (ou plutot une minoration de $\delta(H_{\lambda,i})$). Soit donc $h\in\{0,\ldots, t\}$ maximal tel que $s_h\geq 1$ (on pose $h=0$ si $s_i=0$ pour tout $i$). On a donc 
\[s_1\geq \ldots\geq s_h=1> 0=s_{h+1}=\ldots=s_{t}\ \text{ et }\ P_{r_1,s_1}\subset\ldots P_{r_h,s_h}\subset P_{r_{h+1,0}}\subset\ldots\subset P_{r_{t}}.\]
\noindent Posons
\[\delta_1=\ldots=\delta_h=1\ \text{ et }\ \delta_{h+1}=\ldots=\delta_{t}=0.\]
\noindent Posons $m^{t+1}=0$. On voit (il s'agit d'une somme t\'el\'escopique) que 
\[m^{t+1-h}=m^{t+1-h}-m^{t+1}=\sum_{j=1}^{t}(m^{t+1-j}-m^{t+1-(j-1)})\delta_j.\]
\noindent Or $P_{r_{h},s_h}$ (avec $s_h\geq 1$) correspond au groupe $G_{h,(\lambda,i)}$ lui m\^eme associ\'e \`a l'ensemble $I_{t+1-h}$. Il correspond donc \`a un morceau de $H_{\lambda,i}$ sur lequel on voit qu'il existe $P,Q$ d'ordre $\ell^{m^{t+1-h}}$ tel que l'accouplement de Weil de $\ell^{m^{t-h}}P$ et $\ell^{m^{t-h}}Q$ est une racine primitive $\ell$-eme de $1$. Ceci se traduit en disant que 
\[\delta(H_{\lambda,i})\geq\phi(\ell^{m^{t+1-h}})\que\ell^{m^{t+1-h}},\]
\noindent ceci restant valable pour $h=0$. Nous obtenons ainsi la minoration
\[[K(H_{\lambda,i}):K]\gg\ell^{\sum_{j=1}^{t_{\lambda,i}}(m^{t_{\lambda,i}+1-j}-m^{t_{\lambda,i}+1-(j-1)})(\delta_j+f(\lambda)\codim P_{r_j,s_j})}.\]
\noindent De plus, pour tout entier $k\in\{1,\ldots t\}$,
\[r_{t+1-k}+s_{t+1-k}=|I_k|=\sum_{j=1}^ka_j.\]

\subsection{Invariant $\gamma(A)$ pour $H\subset A[\ell^{\infty}]$}

\noindent Nous sommes ici dans la situation pr\'esent\'ee au d\'ebut de cette section avec $H=\prod_{(\lambda,i)\in\I_{\ell}}X_{\lambda,i}$. Avec les notations introduites dans le cas d'un $H_{\lambda,i}$ (i.e. au paragraphe pr\'ec\'edent \ref{m1}), on peut, pour tout $(\lambda,i)\in\I_{\ell}$, \'ecrire
\[H_{\lambda,i}=\prod_{j=1}^{2h_i}\left(\Z/\ell^{m_j}\Z\right)^{f(\lambda)}=\prod_{j=1}^{t_{\lambda,i}}\left(\Z/\ell^{m^j}\Z\right)^{a_jf(\lambda)},\]
\noindent o\`u $(m^j)_{j\geq 1}$ est une suite strictement d\'ecroissante (le couple $(\lambda,i)$ \'etant fix\'e). 

\medskip

\noindent Nous allons utiliser un r\'esultat galoisien que nous avons d\'emontr\'e dans \cite{hindry-ratazzi1}. Dans le th\'eor\`eme 6.6 de \cite{hindry-ratazzi1} nous donnons une preuve pour $A=\prod_i A_i$ et avec $M_i=T_{\ell}(A_i)$ (cf. notations ci-dessous). En fait la m\^eme preuve reprise mot pour mot donne~: 

\begin{prop}\label{refauprodce}    Soient $T_{\ell}(A)=\oplus_{j\in J}M_j^{\alpha_j}$ une d\'ecomposition de $\Z_{\ell}$-modules galoisiens v\'erifiant les deux propri\'et\'es suivantes (o\`u $M_j[\ell^{\infty}]$ d\'esigne $\cup_n M_j/\ell^nM_j$)~: 
\begin{enumerate}
\item pour tout $j\in J$ et tout groupe fini $H_j\subset M_j[\ell^{\infty}]$, il existe $w_j=w_j(H_j)$ tels qu'on a, uniform\'ement en $(\ell,H_j)$,
\[K(H_j)\cap K(\mu_{\ell^{\infty}})\que K(\mu_{\ell^{w_j}})\ ; \]
\item Uniform\'ement en $\ell$, on a l'identit\'e~:
\[\Gal(K(A[\ell^{\infty}])/K(\mu_{\ell^{\infty}}))\que\prod_{j\in J}\Gal(K(M_j[\ell^{\infty}]))/K(\mu_{\ell^{\infty}}))\]
\end{enumerate}
\noindent Alors si $w:=\max w_j $, pour tout groupe fini $H=\prod_jH_j\subset A[\ell^{\infty}]$, uniform\'ement en $(\ell,H)$ on a, $K(H)\cap K(\mu_{\ell^{\infty}})\que K(\mu_{\ell^w})$ et
\[[K(H):K(\mu_{\ell^{w}})]\que\prod_{j\in J}[K(H_j):K(\mu_{\ell^{w_j}})].\]
\end{prop}

\noindent Nous allons appliquer ceci avec l'ensemble $J=\I_{\ell}$, et pour $j=(\lambda,i)\in\I_{\ell}$, avec $M_j=T_{\lambda,i}(A_i)$, ainsi que $\alpha_j=n_i$ si $A_i$ est de type I et $\alpha_j=2n_i$ si $A_i$ est de type II. Enfin nous l'utiliserons avec $H_j=H_{\lambda,i}$.

\medskip

\noindent Par le lemme \ref{calculH1a} on a 
\[[K(H_{\lambda,i}):K]\que\delta(H_{\lambda,i})(\Hdg(A_i)(\Z_{\ell}):G(H_{\lambda,i})).\]
\noindent Or on sait que dans notre situation on a uniform\'ement en $(\ell,H)$ :
\[\Gal(K(A_i[\ell^{\infty}])/K(\mu_{\ell^{\infty}}]))\que \prod_{(\lambda,i)\in\I_{\ell}}\Gal(K(T_{\lambda,i}[\lambda^{\infty}])/K(\mu_{\ell^{\infty}})).\]
\noindent On peut appliquer la proposition \ref{refauprodce} et on obtient, uniform\'ement en $(\ell,H)$, 
\[[K(H):K]\que\delta(H)\prod_{(\lambda,i)\in\I_{\ell}}(\Sp_{2h_i}(\O_{\lambda,i}):G(H_{\lambda,i})).\]

\newcommand{\cd}{\textnormal{cd}}

\noindent Notons $cd_{(\lambda,i)}^k$ la codimension  du groupe alg\'ebrique $G_{k,(\lambda,i)}$. Dans la situation d'un $H_{\lambda,i}$ fix\'e nous avions introduit au paragraphe pr\'ec\'edent des notations 
\[m_j \text{ et } a_j, \ 1\leq j\leq 2h_i,\ \text{ ainsi que }\ m^r,\ 1\leq r\leq t_{\lambda,i}.\]
\noindent Afin de rendre claire les diverses d\'ependances nous utiliserons ci-dessous les notations un peu plus lourdes suivantes en lieu et place des pr\'ec\'edentes :
\[m_j(\lambda,i) \text{ et } a_j(\lambda,i), \ 1\leq j\leq 2h_i,\ \text{ ainsi que }\ m^r_{\lambda,i},\ 1\leq r\leq t_{\lambda,i}.\]
\noindent Les calculs effectu\'es dans le cas d'un $H_{\lambda,i}$ nous donnent, uniform\'ement en $(\ell,H)$,
\[(\Hdg(A)(\Z_{\ell}):G(H))\que\text{exp}\left(\sum_{(\lambda,i)\in \I_{\ell}}\sum_{k=1}^{t_{\lambda,i}}f(\lambda)\cd_{(\lambda,i)}^k\left(m^{t_{\lambda,i}+1-k}_{\lambda,i}-m^{t_{\lambda,i}+1-(k-1)}_{\lambda,i}\right)\log\ell\right).\]
\noindent De plus, il existe un $(\lambda_1,i_1)$ tel que $\delta(H)=\delta(H_{\lambda_1,i_1})$. Quitte \`a renum\'eroter on peut supposer que $i_1=1$ On note alors $(\delta(\lambda_1)_j)_j$ la suite de $0$ et de $1$ relative \`a $\delta(H_{\lambda_1,1})$ d\'efinie au paragraphe pr\'ec\'edent. On a, uniform\'ement en $(\ell,H)$,
\[\delta(H)\gg \text{exp}\left(\sum_{j=1}^{t_{\lambda_1,1}}\left(m^{t_{\lambda_1,1}+1-j}_{\lambda_1,1}-m^{t_{\lambda_1,1}+1-(j-1)}_{\lambda_1,1}\right)\delta(\lambda_1)_j\log \ell\right).\]
\noindent On pose par ailleurs $\delta(\lambda)_j=0$ pour tout $j$ si $(\lambda,i)\not=(\lambda_1,1)$. Avec ces notations, on trouve en suivant les calculs du cas d'un $H_{\lambda,i}$, la minoration suivante (au sens $\gg$,  uniform\'ement en $(\ell,H)$) pour $[K(H):K]$~:
\[\text{exp}\left(\sum_{(\lambda,i)\in \I_{\ell}}\sum_{j=1}^{t_{\lambda,i}}m^j_{\lambda,i}\left[(\delta(\lambda)_{t_{\lambda,i}+1-j}-\delta(\lambda)_{t_{\lambda,i}+1-(j-1)})+f(\lambda)(\cd_{(\lambda,i)}^{t_{\lambda,i}+1-j}-\cd_{(\lambda,i)}^{t_{\lambda,i}+1-(j-1)})\right]\log\ell\right),\]
\noindent et
\[|H|=\text{exp}\left(\sum_{(\lambda,i)\in \I_{\ell}}n_id_i\sum_{j=1}^{t_{\lambda,i}}m^j_{\lambda,i}f(\lambda)a_j(\lambda,i)\log\ell\right),\]
\noindent o\`u $d_{i}$ vaut 1 (respectivement 2) si $A_i$ est de type I (respectivement de type II) et o\`u l'on rappelle que $A=\prod_{i=1}^dA_i^{n_i}.$

\medskip

\noindent Notons 
\[b_{\lambda,i}^j:=(\delta(\lambda)_{t_{\lambda,i}+1-j}-\delta(\lambda)_{t_{\lambda,i}+1-(j-1)})+f(\lambda)(\cd_{\lambda,i}^{t_{\lambda,i}+1-j}-\cd_{\lambda,i}^{t_{\lambda,i}+1-(j-1)}),\]
\noindent et posons de plus 
\[a_{\lambda,i}^j:=n_id_ia_j(\lambda,i).\]
\noindent Avec ces notations, on aura donc, uniform\'ement en $(\ell,H)$, l'in\'egalit\'e $|H|\ll [K(H):K]^{\gamma}$ si 
\[\gamma\geq\max\frac{\sum_{(\lambda,i)\in\I_{\ell}}\sum_{j=1}^{t_{\lambda,i}}m^j_{\lambda,i}f(\lambda)a_{\lambda,i}^j}{\sum_{(\lambda,i)\in\I_{\ell}}\sum_{j=1}^{t_{\lambda,i}}m^j_{\lambda,i}b_{\lambda,i}^j},\]
\noindent le max \'etant pris sur les $m^1_{\lambda,i}\geq\ldots\geq m^{t_{\lambda,i}}_{\lambda,i}$ pour $(\lambda,i)\in\I_{\ell}$. 

\medskip

\noindent Ainsi, en invoquant le lemme combinatoire \ref{combielem2} et en suivant les notations et calculs du cas d'un $H_{\lambda,i}$, on voit que l'in\'egalit\'e $|H|\ll [K(H):K]^{\gamma}$ est vraie uniform\'ement en $(\ell,H)$, si 
\[\gamma\geq \max \frac{\sum_{(\lambda,i)\in\I_{\ell}} n_id_if(\lambda)(r(\lambda,i)_{t_{\lambda,i}+1-h_{\lambda,i}}+s(\lambda,i)_{t_{\lambda,i}+1-h_{\lambda,i}})}{\delta(\lambda_1)_{t_{\lambda_1,1}+1-h_{\lambda_1,1}}+\sum_{(\lambda,i)\in\I_{\ell}}f(\lambda)\codim P_{r(\lambda,i)_{t_{\lambda,i}+1-h_{\lambda,i}}, s(\lambda,i)_{t_{\lambda,i}+1-h_{\lambda,i}}}}.\]    
\noindent Ce dernier max se r\'e\'ecrit sous la forme
\[\max_{\underset{0\leq s_{\lambda,i}\leq r_{\lambda,i}\leq h_i}{1\leq r_{\lambda,i}}} \frac{\sum_{(\lambda,i)\in\I_{\ell}}n_id_if(\lambda)(r_{\lambda,i}+s_{\lambda,i})}{\delta+\sum_{(\lambda,i)\in\I_{\ell}}f(\lambda) (r_{\lambda,i}+s_{\lambda,i})(2h_i-\frac{1}{2}(r_{\lambda,i}+s_{\lambda,i}-1)},\] 
\noindent et o\`u en reprenant la d\'efinition de $\delta(\lambda_1)_{t_{\lambda_1,1}+1-h_{\lambda_1,1}}$ on voit que $\delta=0$ si tout les $s_{\lambda,i}$ sont nuls et $\delta=1$ si l'un des $s_{\lambda,i}$ est non nul.

\medskip

\noindent Il y a en fait deux \'evaluations \`a faire selon que $\delta=1$ ou que $\delta=0$.

\medskip

\noindent 1. Si $\delta=1$ alors le max \`a \'evaluer se r\'e\'ecrit naturellement sous la forme suivante :
\[\max_{1\leq x_{\lambda,i}\leq 2h_i} \frac{\sum_{(\lambda,i)\in\I_{\ell}}m_id_if(\lambda)x_{\lambda,i}}{1+\sum_{(\lambda,i)\in\I_{\ell}}f(\lambda) x_{\lambda,i}(2h_i-\frac{1}{2}(x_{\lambda,i}-1))}.\] 
\noindent Posons 
\[\rho_1(\underline{x}):=\frac{\sum_{(\lambda,i)\in\I_{\ell}}m_id_if(\lambda)x_{\lambda,i}}{1+\sum_{(\lambda,i)\in\I_{\ell}}f(\lambda) x_{\lambda,i}(2h_i-\frac{1}{2}(x_{\lambda,i}-1))}.\]

\noindent On rappelle que l'on veut comparer $\rho_1(\underline{x})$ avec la quantit\'e
\[\alpha(A):=\max_{I\subset\{1,\ldots,r\}}\frac{2\sum_{i\in I}m_ig_i}{1+\sum_{i\in I}2e_ih_i^2+e_ih_i}.\]

\begin{lemme}pour tout $i\in\{1,\ldots,r\},$ on a $\alpha(A)\geq \frac{m_id_i}{h_i+1}$.
\end{lemme}
\demo C'est un calcul imm\'ediat.\hfill$\Box$

\begin{lemme} \label{comb} On a 
\[\max_{1\leq x_{\lambda,i}\leq 2h_i}\rho_1(\underline{x})\leq \alpha(A).\]
\end{lemme}
\demo L'in\'egalit\'e $\rho_1(\underline{x})\leq \alpha(A)$ se r\'e\'ecrit
\[\sum_{(\lambda,i)\in\I_{\ell}}f(\lambda)\left[x_{\lambda,i}^2-\left(4h_i+1-\frac{2m_id_i}{\alpha(A)}\right)x_{\lambda,i}\right]\leq 2.\]
\noindent Par le lemme pr\'ec\'edent, on voit dans la somme dans le membre de gauche de l'in\'egalit\'e, que les indices tels que $x_{\lambda,i}\leq 2h_i-1$ contribuent via un terme n\'egatif \`a la somme. Autrement dit, la valeur $\rho_1(\underline{x})-\alpha(A)$ est maximale quand pour tout les indices on a $x_{\lambda,i}=2h_i$. Mais dans ce cas, en utilisant que 
\[\sum_{\lambda \text{ place de }\End(A_i)\otimes \Z_{\ell}}f(\lambda)=e_i,\]
\noindent on a:
\begin{eqnarray*}
\rho_1(2h_i,\ldots,2h_i)-\alpha(A)\leq 0	&	\iff	&	\sum_{(\lambda,i)\in\I_{\ell}}f(\lambda)\left[4h_i^2-\left(4h_i+1-\frac{2m_id_i}{\alpha(A)}\right)2h_i\right]\leq 2\\
																					& \iff	& \sum_{i=1}^r\left(-4h_i^2-2h_i+\frac{4m_ih_id_i}{\alpha(A)}\right)\sum_{\lambda }f(\lambda)\leq 2.\\
																					& \iff	&	\sum_{i=1}^r\left(-4h_i^2e_i-2h_ie_i+\frac{4m_ih_ie_id_i}{\alpha(A)}\right)\leq 2.\\
																					& \iff	&	\frac{2\dim A}{\dim\MT(A)}\leq \alpha(A).
\end{eqnarray*}
\noindent La derni\`ere assertion de la s\'erie d'\'equivalences est vraie, ce qui conclut.\hfill$\Box$

\medskip

\noindent 2. Si $\delta=0$ alors dans ce cas un calcul du m\^eme type (plus facile) permet \'egalement de conclure.

\section{Petites valeurs de $\ell$ exceptionnelles}\label{petitl}
\noindent Dans ce paragraphe nous indiquons quelles modifications apporter pour les valeurs exceptionnelles de $\ell$ (en nombre fini). On se place dans le cadre d'une vari\'et\'e ab\'elienne g\'eom\'etriquement simple $A/K$ de type I ou II et pleinement de type Lefschetz, telle que $\End_K(A)=\End_{\bar K}(A)$ et telle que $\End_K(A)$ est un ordre maximal de $D:=\End_K(A)\otimes_{\Z}\Q$. Nous notons $E$ le centre de $D$ et nous notons enfin $\phi : A\rightarrow A^{\vee}$ une polarisation fix\'ee avec $A$ (les diverses constantes intervenant d\'ependant de $A$, d\'ependent aussi du degr\'e de cette polarisation).

\medskip

\noindent Notons que, dans le cas o\`u nous nous sommes plac\'es (pleinement de type Lefschetz), on sait que la conjecture de Mumford-Tate est vraie, donc que l'on a l'inclusion suivante qui est une \'egalit\'e \`a indice fini pr\`es (d\'ependant \'eventuellement de $\ell$ mais peu importe ici car on travaille uniquement avec un nombre fini de valeurs de $\ell$)~:

\[\rho_{\lambda^{\infty}}(G_K)\subset\{x\in\GSp_{2h}(\O_{\lambda})\ |\ \mult(x)\in\Z_{\ell}^{\times}\}.\]

\noindent Nous indiquons dans ce qui suit les petites modifications \`a faire pour pouvoir traiter les $\ell$ qui sont dans l'ensemble fini exceptionnel $S_{\text{ex}}(A)$ introduit dans la d\'efinition 3.4.

\subsection{Si $\ell$ est ramifi\'e dans $\O_E$}

\noindent On suppose dans ce paragraphe que $\ell$ ne divise pas $\deg(\phi)$ et, dans le cas de type II, que $\ell$ est tel que l'alg\`ebre de quaternions $D$ est d\'ecompos\'ee en $\lambda|\ell$. On suppose par contre que $\ell$ est ramifi\'e dans $\O_E~: \ell\O_E=\prod_{\lambda |\ell}\lambda^{e(\lambda)}$.

\medskip

\noindent Nous notons toujours $f(\lambda)$ le degr\'e du corps r\'esiduel en la place $\lambda$. Le lemme \ref{pairing} produit l'accouplement $\phi_{\ell^{\infty}}^{\star}$ sur $T_{\ell}(A)\times T_{\ell}(A)$ \`a valeurs dans $\O_{\ell}^{\star}$. 

\medskip

\noindent \textbf{Hypoth\`ese :~} Supposons pour l'instant pour simplifier que $\phi_{\ell^{\infty}}^{\star}$ est en fait \`a valeurs dans $\O_{\ell}$. Nous verrons plus loin comment faire dans le cas g\'en\'eral. 

\medskip

\noindent Rappelons que dans cette situation on a la d\'ecomposition 
\begin{equation*}
T_{\ell}(A)=\left\{\begin{matrix}\prod_{\lambda|\ell}T_{\lambda}(A)& (\text{Type I})\cr \prod_{\lambda|\ell}T_{\lambda}(A)\oplus T_{\lambda}(A) & (\text{Type II})\cr\end{matrix}\right.
\end{equation*}

\medskip

\noindent Par r\'eduction modulo $\lambda^n$, on obtient alors pour tout entier $+\infty\geq n\geq 1$, comme dans le cas non-ramifi\'e,
\[\phi_{\lambda^n} : T_{\lambda}(A)/\lambda^nT_{\lambda}(A)\times T_{\lambda}(A)/\lambda^n T_{\lambda}(A)\rightarrow O_{\lambda}/\lambda^n(1).\]
\noindent Notons par ailleurs que $\ell O_{\lambda}=\lambda^{e(\lambda)}$, donc par r\'eduction modulo $\ell^n$ on a 
\[ A[\ell^n]=T_{\ell}(A)/\ell^nT_{\ell}(A)=\left\{\begin{matrix}\prod_{\lambda|\ell}T_{\lambda}[\lambda^{e(\lambda)n}]& (\text{Type I})\cr \prod_{\lambda|\ell}T_{\lambda}[\lambda^{e(\lambda)n}]\oplus T_{\lambda}[\lambda^{e(\lambda)n}]& (\text{Type II})\cr\end{matrix}\right. ,\]
\noindent et 
\[\phi_{\lambda^{e(\lambda)n}} : T_{\lambda}[\lambda^{e(\lambda)n}]\times T_{\lambda}[\lambda^{e(\lambda)n}]\rightarrow \O_{\lambda}/\ell^n \O_{\lambda}(1).\]
\noindent De plus on v\'erifie que 
\[\phi_{\lambda^{e(\lambda)n}}(\ell x,\ell y)=\phi_{\lambda^{e(\lambda)(n+1)}}(x,y)^{\ell}\]
\noindent et on voit que l'action de Galois sur les accouplements $\phi_{\lambda^{e(\lambda)n}}$ se fait via le caract\`ere cyclotomique $\chi_{\ell^n}$. Par ailleurs, le $\Z_{\ell}$-module $\O_{\lambda}$ \'etant libre de rang $e(\lambda)f(\lambda)$, on a
\begin{equation}\label{eq1}
O_{\lambda}/\ell^n\O_{\lambda}\simeq \left(\Z/\ell^n\Z\right)^{e(\lambda)f(\lambda)}\text{ et } T_{\lambda}[\lambda^{e(\lambda)n}]\simeq \left(O_{\lambda}/\ell^n\O_{\lambda}\right)^{2h}\simeq \left(\Z/\ell^n\Z\right)^{2he(\lambda)f(\lambda)}.
\end{equation}

\noindent Finalement en travaillant en lieu de place de $\phi_{\lambda^n}$ avec $\phi_{\lambda^{e(\lambda)n}}$, en d\'ecomposant $H$ selon les $H_{\lambda}\subset T_{\lambda}[\lambda^{e(\lambda)n}]\subset A[\ell^n]$, on peut reprendre tout ce qui a \'et\'e fait dans le cas non-ramifi\'e. Notamment la propri\'et\'e $\mu$ pour les $H_{\lambda}$ est toujours v\'erifi\'ee avec la modification \'evidente suivante dans la d\'efinition de $m_1(H_{\lambda})$ (et de m\^eme pour $m(H_{\lambda})$) : on pose 
\[m_1(H_{\lambda}):= \max\left\{m\ |\ \exists P,Q\in H_{\lambda} \text{ d'ordre }\ell^k \text{ tels que } \phi_{\lambda^{e(\lambda)k}}(P,Q)\text{ est d'ordre }\ell^m\right\}.\]
\noindent Pour prouver la propri\'et\'e $\mu$ dans ce cadre on peut toujours utiliser le paragraphe sur les modules isotropes sur $(\O_{\lambda},\F_{\lambda})$, puisque les accouplements $\lambda$-adiques $\phi_{\lambda^r}$ sont construits modulo $\lambda^{r}$ pour tout entier $r\geq 1$ et non pas uniquement modulo $\ell^{r}$ (cette r\'eduction modulo $\ell^{r}$ \'etant celle utilis\'ee pour d\'efinir les $H_{\lambda}$).

\medskip

\noindent En reprenant mot pour mot les calculs combinatoires d\'ej\`a effectu\'es, on voit avec (\ref{eq1}) que les calculs restent inchang\'es sous r\'eserve de remplacer partout $f(\lambda)$ par $f(\lambda)e(\lambda)$. La contrainte $\sum f(\lambda)=[E:\Q]$ \'etant remplac\'ee par la contrainte $\sum e(\lambda)f(\lambda)=[E:\Q]$, on voit que la valeur $\gamma$ reste la m\^eme dans ce cadre, ce que l'on voulait prouver.

\medskip

\noindent \textbf{Suppression de l'hypoth\`ese :} on ne suppose d\'esormais plus que $\phi_{\ell^{\infty}}^{\star}$ est \`a valeurs dans $\O_{\ell}$. 

\medskip 

\noindent On sait que $\O_{\lambda}^{\star}$ est un id\'eal fractionnaire de $E_{\lambda}$, donc de la forme $\pi_{\lambda}^{-m_{\lambda}}\O_{\lambda}$ pour un certain entier $m_{\lambda}$ (avec $\pi_{\lambda}$ une uniformisante). On choisit 
\[m_0:=\text{pgcd}\left(m_{\lambda}\ |\ \lambda|\ell,\ \ell\text{ ramifi\'e dans }\O_{E}\right).\]
\noindent On a ainsi $\ell^{m_0}\O_{\ell}^{\star}\subset \O_{\ell}$ pour tout les $\ell$ que l'on consid\`ere. On fait alors les modifications suivantes : 
\begin{enumerate}
\item On remplace au d\'epart le module $T_{\ell}(A)$ par $T_{\ell}':=\ell^{m_0}T_{\ell}(A)$. Dans ce cas, le lemme \ref{pairing} produit l'acccouplement $\phi_{\ell^{\infty}}^{\star}$ sur $T_{\ell}'\times T_{\ell}'$ \`a valeurs dans $\O_{\ell}$ (c'est pour arriver dans $\O_{\ell}$ que l'on a remplac\'e $T_{\ell}(A)$ par $T_{\ell}'$). 
\item On travaille avec $H'=\ell^{m_0}H$ en lieu et place de $H$. La raison de cette modification est que le groupe $H$ n'est a priori pas contenu dans la r\'eduction modulo $\ell^n$ de $T_{\ell}'$. Par contre $H'$ l'est.
\end{enumerate}
\noindent \`A la d\'eperdition pr\`es d'indice en $\ell^{2hm_0}$ pr\`es, on peut reprendre la preuve d\'ej\`a effectu\'ee et on obtient, uniform\'ement en $H$ : 
\[|H'|\ll[K(H') :K]^{\gamma(A)}\leq [K(H):K]^{\gamma(A)}.\]
\noindent De plus, on a visiblement $|H|\leq |H'|\cdot|A[\ell^{m_0}]|$. Les premiers probl\'ematiques $\ell$ \'etant en nombre born\'es, il en est de m\^eme pour le cardinal des divers $A[\ell^{m_0}]$ et on voit donc que l'on obtient ainsi, uniform\'ement en $H$,
\[|H|\ll [K(H):K]^{\gamma(A)}.\]

\subsection{Si $\ell$ divise le degr\'e de la polarisation}
\noindent On suppose d\'esormais que $\ell$ est un premier quelconque divisant $\deg(\phi)$ et, dans le cas de type II, que $\ell$ est tel que l'alg\`ebre de quaternions $D$ est d\'ecompos\'ee en $\lambda|\ell$. Notons $m_0$ l'entier maximal tel que $\ell^{m_0}$ divise $\deg(\phi)$. Dans ce cas, toute les constructions faites jusqu'\`a pr\'esent continuent encore \`a s'appliquer \`a condition de faire au d\'epart les modifications suivantes :
\begin{enumerate}
\item On travaille avec $T_{\ell}'(A):=\ell^{m_0}T_{\ell}(A)$ en lieu et place de $T_{\ell}(A)$.
\item On travaille avec l'accouplement $\phi_{\ell}^{\phi} : T_{\ell}'(A)\times T_{\ell}'(A)\rightarrow  T_{\ell}'(A)\times T_{\ell}'(A^{\vee})\rightarrow \Z_{\ell}(1)$ d\'efini par $(x,y)\mapsto \phi_{\ell}(x,\phi(y))$, en lieu et place de l'accouplement $\phi_{\ell}$. Ce choix ainsi que le point pr\'ec\'edent sont faits de sorte \`a avoir un accouplement bilin\'eaire alern\'e sur $T_{\ell}'(A)$, non d\'eg\'en\'er\'e modulo $\ell^n$ pour tout $n\geq 1$.
\item On travaille avec $H'=\ell^{m_0}H$ en lieu et place de $H$. La raison de cette modification est que le groupe $H$ n'est a priori par contenu dans la r\'eduction modulo $\ell^n$ de $T_{\ell}'(A)$. Par contre $H'$ l'est.
\end{enumerate}
\noindent Avec ces modifications on peut reprendre la preuve d\'ej\`a effectu\'ee et on obtient, uniform\'ement en $H$~: 
\[|H'|\ll[K(H') :K]^{\gamma(A)}\leq [K(H):K]^{\gamma(A)}.\]
\noindent De plus, on a visiblement $|H|\leq |H'|+|A[\ell^{m_0}]|$. Les premiers $\ell$ \'etant en nombre fini, le le cardinal des divers $A[\ell^{m_0}]$ est donc bornée et l'on obtient ainsi uniform\'ement en $H$,
\[|H|\ll [K(H):K]^{\gamma(A)}.\]

\subsection{Petites valeurs exceptionnelles pour le type II}
\noindent Dans le cas d'une vari\'et\'e ab\'elienne de type II, il y a encore un nombre fini de valeurs $\ell$ exceptionnelles \`a traiter : les premiers $\ell$ tels que l'alg\`ebre $D$ est non-d\'ecompos\'ee en $\lambda|\ell$. Ce cas des $\ell$ ramifi\'es dans l'alg\`ebre de quaternions doit \^etre trait\'e avec une l\'eg\`ere modification : la d\'ecomposition 
 $V_{\ell}(A)=\prod_{\lambda}\left(W_{\lambda}(A)\oplus W_{\lambda}(A)\right)$ n'existant que apr\`es tensorisation par une extension quadratique (voir Proposition \ref{2copies}). Avec les notations du d\'ebut du paragraphe 7 : au lieu d'avoir sur $\O_{\lambda,i}$ la d\'ecomposition $X_{\lambda,i}=H_{\lambda,i}\oplus H_{\lambda,i}$, on a, en travaillant sur une extension quadratique de $\O_{\lambda,i}$, la d\'ecomposition $X_{\lambda,i}=H_{\lambda,i}\oplus \bar H_{\lambda,i}$ o\`u $\bar H_{\lambda,i}$ est conjugu\'e \`a $H_{\lambda,i}$. Ainsi au lieu de comparer le cardinal de $H_{\lambda,i}$ avec le degr\'e de l'extension $[K(H_{\lambda,i}):K]$ on reprend la m\^eme preuve en travaillant directement avec $X_{\lambda,i}$, comparant le cardinal de $X_{\lambda,i}$ et le degr\'e de $[K(X_{\lambda,i}):K]$.

\section{Ordre d'un point et degr\'e de l'extension qu'il engendre}

\noindent Nous donnons dans ce paragraphe la preuve du th\'eor\`eme \ref{unpoint}. Nous pouvons pour cela supposer (et nous le faisons) que tous les $\bar K$-endomorphismes de $A$ sont d\'efinis sur $K$.

\medskip

\noindent Soit $P$ un point de torsion et $H_P$ le $\End(A)$-module engendr\'e par $P$, on a clairement $K({P})=K(H_P)$. En remarquant que $\codim P_{1,0}=2h$ (dans $\Sp_{2h}$), les arguments des paragraphes \ref{p6} et \ref{p7} pr\'ec\'edents montrent qu'un point  $P_{\lambda}$ d'ordre $\ell^n$ dans $T_{\lambda}[\lambda^{\infty}]$ engendre (uniform\'ement en ($\ell,P$)) une extension de degr\'e
\[ [K(P_{\lambda}):K]\gg  \ell^{2hn}.\]
\noindent Si ensuite $P=\sum_{\lambda\,|\,\ell}P_{\lambda}$ avec $P_{\lambda}$ point de $T_{\lambda}[\lambda^{\infty}]$  et d'ordre $\ell^{n_{\lambda}}$, de sorte que $P$ est d'ordre $\ell^n$ avec $n=\max n_{\lambda}$, alors, uniform\'ement en ($\ell,P$) on a~:
\[ [K({P}):K]\gg\ell^{2h\sum_{\lambda}n_{\lambda}}\geq \ell^{2hn}.\]
\noindent Enfin si $P$ est d'ordre $m$ quelconque avec $m=\prod_{i=1}^r\ell_i^{n_i}$, on peut \'ecrire $P=\sum_{i=1}^rP_i$ avec $P_i$ d'ordre $\ell_i^{n_i}$. L'ind\'ependance des repr\'esentations $\ell$-adiques (Cf proposition \ref{indep}) permet d'\'ecrire, uniform\'ement en ($\ell,P$), 
\[ [K({P}):K]=[K({P_1},\dots,{P_r}):K]\gg \prod_{i=1}^r[K({P_i}):K]\geq \prod_{i=1}^rc_1\ell^{2hn_i}=c_1^{\omega(m)}m^{2h}.\]

\section{Appendice : compl\'ements autour de la conjecture de Mumford-Tate\label{aboutMT}}

\subsection{Indice de l'image de Galois dans le groupe de Mumford-Tate}

La conjecture de Mumford-Tate dit que l'inclusion $G^0_{\ell^{\infty}}\subset\MT\otimes\Q_{\ell}$ est une \'egalit\'e, ou encore que, quitte a avoir effectu\'e une extension finie du corps de base $K$, l'image de la repr\'esentation $\ell$-adique galoisienne $\rho(G_K)$ est contenue et ouverte dans $\MT(\Q_{\ell})$, ou encore comme 
$\GL(T_{\ell}(A))\cong \GL_{2g}(\Z_{\ell})$ est compact, la conjecture \'equivaut \`a dire que, quitte a avoir effectu\'e une extension finie du corps de base $K$, $\rho(G_K)$ est d'indice fini dans $\MT(\Z_{\ell})$. Une forme l\'eg\`erement plus forte, sugg\'er\'ee par Serre, affirme que cet indice est born\'e ind\'ependamment de $\ell$. Clarifions tout d'abord ce point en montrant que la conjecture de Mumford-Tate entra\^{\i}ne la forme ``forte".

\begin{theo}\label{indicefini} Si $A/K$ v\'erifie la conjecture de Mumford-Tate alors l'indice de $\rho(G_K)$ dans $\MT(A)(\Z_{\ell})$ est born\'e ind\'epen\-damment de $\ell$.
\end{theo}
\noindent La preuve consiste \`a r\'eunir un r\'esultat de Serre \cite{serrekr86} (resp. de Wintenberger \cite{W}) concernant la partie torique centrale (resp. la partie semi-simple) des groupes $\ell$-adiques et des groupes de Mumford-Tate. Plus pr\'ecis\'ement, notons $S=S_A$ le groupe d\'eriv\'e du groupe de Mumford-Tate de $A$ ou, ce qui revient au m\^eme, du groupe de Hodge et notons $C$ la composante neutre du centre du groupe $\MT(A)$; ce sont des $\Q$-groupes alg\'ebriques. Notons similairement $S_{\ell}=S_{\ell,A}$ le groupe d\'eriv\'e de $G_{\ell,A}$ ou, ce qui revient au m\^eme, du groupe $H_{\ell,A}$ et notons $C_{\ell}$ la composante neutre du centre du groupe $G_{\ell,A}$; ce sont des $\Q_{\ell}$-groupes alg\'ebriques.
 
\noindent On sait par les travaux de Borovo{\u\i} \cite{boro}, Deligne \cite{deligne} Exp I, 2.9, 2.11, et Pjatecki{\u\i}-{\v{S}}apiro \cite{piat} que
\begin{equation}
C_{\ell}\subset C_{\Q_{\ell}}\qquad {\rm et}\qquad S_{\ell}\subset S_{\Q_{\ell}}.
\end{equation}
\noindent En fait on sait m\^eme que la premi\`ere inclusion est une \'egalit\'e, essentiellement d'apr\`es la th\'eorie ab\'elienne de Serre \cite{mg}, une preuve est d\'etaill\'ee dans \cite{vasiu} et reprise dans \cite{ullmoyaf}. On ne sait pas, en g\'en\'eral si la deuxi\`eme inclusion est une \'egalit\'e, en fait l'\'egalit\'e est  \'equivalente \`a la conjecture de Mumford-Tate. D'apr\`es Faltings, les deux groupes r\'eductifs ont le m\^eme commutant, leur \'egalit\'e est aussi \'equivalente \`a l'\'egalit\'e des {\it rangs} des deux groupes semi-simples.

\medskip
 
\noindent Posons $V_{\ell}(A)=T_{\ell}(A)\otimes_{\Z_{\ell}}\Q_{\ell}$. Les groupes $C$ et $S$ sont des sous-groupes du groupe de Mumford-Tate $\MT=\MT_A$. En voyant $\MT_{\Q_{\ell}}$ comme un sous-groupe alg\'ebrique de $\GL_{V_{\ell}(A)}\cong\GL_{2g,\Q_{\ell}}$, on peut \'etendre ces groupes sur $\Z_{\ell}$ en prenant leur adh\'erence de Zariski dans   $\GL_{T_{\ell}(A)}$. Avec un l\'eger abus de notation nous noterons $C(\Z_{\ell})$ (resp. $S(\Z_{\ell})$) le groupe des $\Z_{\ell}$-points de $C$ (resp. $S$) vu comme groupe alg\'ebrique sur $\Z_{\ell}$. Le m\^eme proc\'ed\'e nous permet d'\'etendre $C_{\ell}$ (resp. $S_{\ell}$) en un groupe sur $\Z_{\ell}$.
 \medskip
 
\noindent Concernant la partie torique centrale, nous savons donc que $\rho_{\ell^{\infty}}(G_K)\cap C(\Z_{\ell})$ est d'indice fini dans $C(\Z_{\ell})$. Le r\'esultat suivant de Serre pr\'ecise ce point et est un des deux points clef pour la preuve du th\'eor\`eme \ref{indicefini}.

\begin{prop}\label{indicab}  \textnormal{\textbf{(Serre) \cite{serrekr86}}} L'indice $\left(C(\Z_{\ell}):C_{\ell}(\Z_{\ell})\cap \rho_{\ell^{\infty}}(G_K)\right)=:c_{\ell}$ est fini born\'e ind\'epen\-damment de $\ell$.
\end{prop}
\noindent  C'est le th\'eor\`eme p.60 de \cite{serrekr86}. Notons que la preuve  de \cite{serrekr86} est r\'edig\'ee   dans le cas $\End(A)=\Z$ et esquiss\'ee dans le cas g\'en\'eral. Pour la commodit\'e du lecteur nous   donnons  ci-dessous une description du centre et de la relation avec l'image de Galois.

\medskip

\noindent Concernant la partie semi-simple $S_{\ell}$ de $G_{\ell}$, suivant Wintenberger \cite{W}, notons $S_{\ell, sc}\rightarrow S_{\ell}$ le rev\^etement universel de $S_{\ell}$ (sur $\Z_{\ell}$) et posons
\[S_{\ell}(R)_u\ \text{ l'image de } S_{\ell, sc}(R)\ \text{ dans }\ S_{\ell}(R),\ \text{ pour }\ R\in\{\Z_{\ell}, \Q_{\ell}, \F_{\ell}\}.\]
\noindent Si $\ell\geq 5$ alors le groupe $S_{\ell}(\F_{\ell})_u$ est le sous-groupe de $S_{\ell}(\F_{\ell})$ engendr\'e par les \'el\'ements unipotents. C'est \'egalement le groupe des commutateurs de $S_{\ell}(\F_{\ell})$. Le point clef que nous utilisons peut s'\'enoncer ainsi.

\begin{prop}\label{wint}\textnormal{\textbf{(Wintenberger) \cite{W}}} L'indice $S_{\ell}(\Z_{\ell})_u$ dans $S_{\ell}(\Z_{\ell})$ est born\'e ind\'ependamment de $\ell$.   Pour tout premier $\ell$ assez grand, le groupe $S_{\ell}(\Z_{\ell})_u$ est contenu dans $G_{\ell}$.
 En particulier, on a les inclusions $S_{\ell}(\Z_{\ell})_u\subset G_{\ell}\cap S_{\ell}(\Z_{\ell})\subset S_{\ell}(\Z_{\ell})$ avec   indices finis, born\'es ind\'ependamment de $\ell$.
\end{prop}

\noindent Cet \'enonc\'e d\'ecoule de l'\'enonc\'e plus pr\'ecis suivant.

\begin{prop}\label{wintenb}\textnormal{\textbf{(Wintenberger) \cite{W}}} On a l'\'egalit\'e $S_{\ell}(\Z_{\ell})_u=S_{\ell}(\Z_{\ell})\cap S_{\ell}(\Q_{\ell})_u$. De plus, si le centre $Z(S_{\ell, sc})$ est de cardinal premier \`a $\ell$ alors $S_{\ell}(\Z_{\ell})_u$ est l'image r\'eciproque de $S_{\ell}(\F_{\ell})_u$ par le morphisme de r\'eduction modulo $\ell$, $\pi_{\ell} : S_{\ell}(\Z_{\ell})\rightarrow S_{\ell}(\F_{\ell}).$
Pour tout premier $\ell$ assez grand, le groupe $S_{\ell}(\Z_{\ell})_u$ est contenu dans $G_{\ell}$. Si $\ell$ est assez grand, alors l'indice $(S_{\ell}(\Z_{\ell}) : S_{\ell}(\Z_{\ell})_u)$ est major\'e par  $c(2\dim A):=\text{ppcm}(n\ |\ n\leq 2\dim A)$.
\end{prop}

\noindent La d\'emonstration du th\'eor\`eme \ref{indicefini} est maintenant imm\'ediate \`a partir des propositions \ref{indicab} et \ref{wint}. Comme $S$ (resp. $C$) est  le groupe d\'eriv\'e (resp. la composante neutre du centre) de $\MT(A)$ on a $\MT(A)=C\cdot S$ et on en tire ais\'ement que $(\MT(A)(\Z_{\ell}):C(\Z_{\ell})\cdot S(\Z_{\ell}))$ est born\'e ind\'ependamment de $\ell$. La proposition \ref{indicab} fournit un sous-groupe $C_1$ de $C(\Z_{\ell})\cap\rho_{\ell^{\infty}}(G_K)$  d'indice fini dans $C(\Z_{\ell})$,   tandis que la proposition \ref{wint}, jointe \`a l'hypoth\`ese 
$S_{\Q_{\ell}}=S_{\ell}$ fournit un sous-groupe $S_1$ de $S(\Z_{\ell})\cap\rho_{\ell^{\infty}}(G_K)$  d'indice fini dans $S(\Z_{\ell})$. On conclut  bien alors que
$(MT(A)(\Z_{\ell}):\rho(G_K))\leq (MT(A)(\Z_{\ell}):C_1\cdot S_1)$ est born\'e ind\'ependamment de $\ell$.

\medskip

\noindent Donnons maintenant la description promise du centre du groupe de Mumford-Tate.

\medskip

\noindent Notons $L=\prod_iL_i$ le centre de $\End^0(A)$; chaque $L_i$ est un corps de nombres et la d\'ecomposition correspond \`a la d\'ecomposition de $A$ \`a isog\'enie pr\`es en composantes isotypiques, i.e. $A\cong \prod_i A_i$ avec $A_i=B_i^{m_i}$ et $B_i$ absolument simple.
On pose aussi $T_L=\prod_i\Res_{L_i/\Q}(\G_{m,L_i})$ et on note $\det_{L_i}:\Res_{L_i/\Q}\left(\GL_{V_i,L_i}\right)\rightarrow\Res_{L_i/\Q}\left(\G_{m,L_i}\right)$ et $\det_L=\prod_i\det_{L_i} $. La restriction de $\det_L$ au tore $T_L$ donne une isog\'enie $\delta:T_L\rightarrow T_L$ qui peut \^etre explicit\'ee comme l'application $x=(x_i)_{i\in I}\mapsto (x_i^{d_i})_{i\in I}$, o\`u $d_i:=2\dim A_i/[L_i:\Q]$. Introduisons une extension auxiliaire $\tilde{L}$ finie, galoisienne sur $\Q$ et contenant les $L_i$; on d\'efinit ensuite la ``norme"  (Cf page 135 de \cite{ichi}):
$$\Psi_i:\Res_{\tilde{L}/\Q}\G_{m,\tilde{L}}\rightarrow\Res_{L_i/\Q}\G_{m,L_i}\quad{\rm et}\quad \Psi=\prod_i\Psi_i:\Res_{\tilde{L}/\Q}\G_{m,\tilde{L}}\rightarrow\prod_i\Res_{L_i/\Q}\G_{m,L_i}.$$

\noindent On a alors une description de la composante neutre du centre du groupe de Mumford-Tate comme le sous-tore de $T_L$ v\'erifiant $\delta({C})={\rm Im}\,\Psi$ (aux notations pr\`es, c'est la proposition 1.2.1 de \cite{ichi}).

\medskip
   
\begin{prop} (Cf. \cite{ichi,serrekr86}) {\rm Le tore $C=C_A$ est le sous-tore de $T_L$ tel que $\delta({C})=\Psi(T_L)$.}  
\end{prop}

\medskip

\noindent Le lien avec les repr\'esentations $\ell$-adiques peut \^etre d\'ecrit ainsi (Cf Serre \cite{serremfv}).

\medskip


 
\noindent Chaque morceau $V_{\ell}(A_i)$ est libre de rang $d_i$ sur $L_i\otimes\Q_{\ell}$; si l'on pose $L_{\ell}:=L\otimes\Q_{\ell}$, alors  $V_{\ell}(A)$ est un  $L_{\ell}$-module et on peut d\'efinir $\det_{L_{\ell}}(V_{\ell}(A))$ qui est libre de rang 1, ce qui fournit une repr\'esentation $\ell$-adique ab\'elienne \`a valeurs dans $T_L$~:
\[\phi_{\ell}:\,G_K\rightarrow L_{\ell}^{\times}=(L\otimes\Q_{\ell})^{\times}=T_L(\Q_{\ell}).\]
\noindent  Par la th\'eorie ab\'elienne de Serre il existe un module $\frak{m}$ et un  homomorphisme de groupes alg\'ebriques associ\'e $S_{\frak{m}}\rightarrow T_L$ induisant la repr\'esentation de la mani\`ere suivante; rappelons (Cf \cite{mg}) que $ S_{\frak{m}}$ est un $\Q$-groupe alg\'ebrique extension du groupe fini $C _{\frak{m}}$ des classes d'id\`eles  modulo $\frak{m}$ par un tore $ T_{\frak{m}}$. Si on note $\epsilon_{\ell}:G_K\rightarrow S_{\frak{m}}(\Q_{\ell})$ la repr\'esentation d\'efinie dans \cite{mg}, alors  $\phi_{\ell}=\phi\circ\epsilon_{\ell}$.
  
\medskip
  
\noindent Le lien avec le centre du groupe de Hodge est que $\Psi(T_L)$ est la composante neutre de $\phi(S_{\frak{m}})$. Ce fait appel\'e ``exercice emb\^etant" dans \cite{serremfv} est \'equivalent \`a l'\'egalit\'e de la composante connexe des centres de $G_{\ell}$ et $\MT(A)_{\Q_{\ell}}$ cit\'ee ci-dessus. 

\subsection{Quelques cas de la conjecture de Mumford-Tate}

\noindent Pour \'enoncer le premier r\'esultat en vue, nous rappelons la notation suivante, Cf. (\ref{defdeS})~: 
\begin{equation}
\Sigma=\left\{g\geq 1\;|\; \exists k\geq 3,\;{\rm impair},\; \exists a\geq 1,\; 2g=(2a)^{k} \text{ ou } 2g={2k\choose k}\right\}.
\end{equation}

\begin{theo}\label{t1} Soit $h$ un entier tel que $h\notin\Sigma$ et soit $A/K$ une vari\'et\'e ab\'elienne de type II telle que le centre de $\End(A)\otimes\Q$ est r\'eduit \`a $\Q$. On suppose que $A$ est de dimension $g=2h$. La conjecture de Mumford-Tate est vraie pour $A$ et le groupe de Mumford-Tate associ\'e \`a $A$ est $\GSp_{2h,\Q}$.
\end{theo}
\demo Il s'agit tout simplement d'appliquer la proposition 4.7 de Pink \cite{pink}. Pr\'ecis\'ement : soit $\ell$ un nombre premier suffisamment grand ; d\'ecomposons le module de Tate $V_{\ell}$ en somme de 2 copies isomorphes $W_{\ell}$ sur $\Q_{\ell}$. Les $W_{\ell}$ donnent des repr\'esentations de dimension $2h$, de type Mumford-Tate forte (puisque c'est le cas par \cite{pink} theorem 5.10 pour les $V_{\ell}$), fid\`ele, symplectique, absolument irr\'eductible du groupe d\'eriv\'e de $G_{\ell}$. On peut donc appliquer la proposition 4.7 de Pink qui nous dit que si $h$ est en dehors de l'ensemble exceptionnel $\Sigma$ alors $G_{\ell}$ est isomorphe \`a $\GSp_{2h,\Q}$.
\hfill$\Box$

\bigskip
Pour \'enoncer le r\'esultat suivant, rappelons que si une vari\'et\'e ab\'elienne $A$ a r\'eduction semi-stable en une place $v$, la composante neutre de la fibre sp\'eciale est une extension d'une vari\'et\'e ab\'elienne $B_0$ par un tore que nous noterons $T_0$; la dimension de ce tore s'appelle la {\it dimension torique}. De plus, si ce tore est non trivial (i.e. s'il y a mauvaise r\'eduction) l'anneau des endomorphismes de $A$ agit fid\`element sur ce tore. En cons\'equence, si $A$ est de type I et $e=[\End^0(A):\Q]$ (resp. de type II et  $4e=[\End^0(A):\Q]$) alors la dimension de $T_0$ est un multiple de $e$ (resp. de $2e$).

\begin{theo}\label{hallbis} Soit $A$ une vari\'et\'e ab\'elienne simple de type I ou II, d\'efinie sur un corps de nombres $K$. On note $e$ le degr\'e sur $\Q$ du centre de $\End^0(A)$ et $d:=1$, si $A$ est de type I (resp. $d:=2$ si $A$ est de type II).  Supposons de plus que:

$(T):\quad$ il existe une place de $K$ o\`u $A$ poss\`ede   r\'eduction semi-stable de dimension torique $de$.

Alors $A$ est pleinement de type Lefschetz.
\end{theo}
Commen\c{c}ons par rappeler quelques r\'esultats de Grothendieck d\'ecrits dans \cite{groth}. On \'etudie la r\'eduction modulo un id\'eal premier donc on peut passer \`a un anneau local que l'on peut d'ailleurs compl\'eter et qu'on notera $\mathcal{O}$; dans le langage de \cite{groth}, le spectre de $\mathcal{O}$ est un {\it trait}. Tout sch\'ema quasi-fini $X/\mathcal{O}$ se d\'ecompose en
$$X=X^{\rm f}\sqcup X',$$
o\`u $X^{\rm f}$ est un sch\'ema fini sur $\mathcal{O}$ et $X'$ est un sch\'ema \'egal \`a sa fibre g\'en\'erique. On choisit un nombre premier $\ell$ assez grand pour qu'il ne divise pas le cardinal du groupe des composantes de la fibre sp\'eciale du mod\`ele de N\'eron sur $\mathcal{O}$ de $A$. En d\'ecomposant ainsi le sch\'ema quasi-fini $A[\ell^n]$ et en prenant la limite on obtient la ``{\it partie fixe}" 
$T_{\ell}(A)^{\rm f}\subset T_{\ell}(A)$.
En consid\'erant la partie torique de la fibre sp\'eciale $A_0^0$ du mod\`ele de N\'eron
$$0\rightarrow T_0\rightarrow A_0^0\rightarrow B_0\rightarrow 0,$$
on obtient la ``{\it partie torique}" 
$T_{\ell}(A)^{\rm t}\subset T_{\ell}(A)^{\rm f}$, qui est de rang $\mu:=\dim T_0$ (i.e. \'egal \`a la dimension torique). Le r\'esultat clef de Grothendieck d\'ecrit ces sous-modules en termes de l'action du groupe d'inertie, not\'e $I$ et s'\'enonce ainsi (le premier point est la {Proposition 2.2.5} de \cite{groth}, le deuxi\`eme le th\'eor\`eme {Th\'eor\`eme 2.4} de \cite{groth}).


\begin{theo} \textnormal{\textbf{(Grothendieck \cite{groth})}} Avec les notations pr\'ec\'edentes, on a 
\begin{enumerate}
\item On a l'\'egalit\'e $T_{\ell}(A)^{\rm f}= T_{\ell}(A)^{I}$.
\item Soit $\check{A}$ la duale de $A$, d\'esignons par $\perp$ l'orthogonal au sens de l'accouplement canonique de Weil $T_{\ell}(A)\times T_{\ell}(\check{A})\rightarrow T_{\ell}(\G_m)$. Alors 
$$T_{\ell}(A)^{\rm t}= T_{\ell}(A)^{\rm f}\cap\left( T_{\ell}(\check{A})^{\rm f}\ \right)^{\perp}=T_{\ell}(A)^{I}\cap \left(T_{\ell}(\check{A})^{I}\right)^{\perp}.$$
\end{enumerate}
\end{theo} 
De plus, si on note $U=V_{\ell}(A)=T_{\ell}(A)\otimes_{\Z_{\ell}}\Q_{\ell}$, $U_1=V_{\ell}(A)^{\rm f}$ et $U_2=V_{\ell}(A)^{\rm t}$, et si l'on identifie $V_{\ell}(A)$ et $V_{\ell}(\check{A})$ via une polarisation fix\'ee, on a donc $U_2=U_1\cap U_1^{\perp}$. L'inertie agit de plus trivialement sur $U_2$ et sur $U/U_2$ (mais non trivialement sur $U$ si on suppose qu'il y a mauvaise r\'eduction).

\medskip

\demo (du th\'eor\`eme \ref{hallbis}). Notons comme pr\'ec\'edemment $E$ le centre de $\End^0(A)$.
Il suffit de montrer que, pour un $\ell$, le groupe $H_{\ell}$ est le produit des $\Sp_{2h,E_{\lambda}}$, quand $\lambda$ parcourt les id\'eaux premiers de $E$ au dessus de $\ell$. On choisit un $\ell$ assez grand et totalement d\'ecompos\'e dans $E/\Q$; d'apr\`es le lemme de rel\`evement (lemme \ref{releveserre}), il suffit de voir que l'image de Galois modulo $\ell$ contient le produit des $\Sp_{2h,\F_{\ell}}$. Chris Hall \cite{hall} montre dans le cas o\`u $\End(A)=\Z$ (ie $e=d=1$) que l'hypoth\`ese $(T)$ entra\^ine ceci. Nous expliquons comment adapter ses arguments au cas plus g\'en\'eral.

Consid\'erons les sous-modules d\'ecrits ci-dessus $U_2\subset U_1\subset V_{\ell}(A)$ et rappelons que $\mu=\rang U_2$ est la dimension torique.
L'inertie op\`ere non trivialement sur $V_{\ell}(A)$ puisqu'il y a mauvaise r\'eduction mais trivialement sur $U_2$ et $V_{\ell}(A)/U_2$. Si $g_I$ d\'esigne un g\'en\'erateur topologique du quotient maximal pro-$\ell$-fini de $I$, sa matrice $\rho_{\ell}(g_I)$ est donc conjugu\'ee \`a une matrice $\begin{pmatrix}  I_{\mu}& B\cr 0&I_{2g-\mu}\cr\end{pmatrix}$.

On peut ``d\'ecouper" les  $U_i$, qui sont des $\Q_{\ell}$-espaces vectoriels, en introduisant $e_{\lambda}$ l'idempotent de $E_{\ell}=E\otimes_{\Q}\Q_{\ell}$ qui projette $E_{\ell}$ sur $E_{\lambda}$ et en notant $U_{i,\lambda}$ l'image $e_{\lambda}U_i$.  On obtient alors la d\'ecomposition cherch\'ee.
\begin{lemme} Les sous-modules galoisiens $U_1$  ``partie fixe" et $U_2$ ``partie torique" sont stables par $\End(A)$ et se d\'ecomposent ainsi.
\begin{enumerate}
\item (Type I) Les $U_i$  se d\'ecomposent    en $U_i=\prod_{\lambda}U_{i,\lambda}$.
\item (Type II) Les $U_i$  se d\'ecomposent    en $U_i=\prod_{\lambda}\left(U_{i,\lambda}\oplus U_{i,\lambda}\right)$.
\end{enumerate}
\end{lemme}
\demo (du lemme). L'action des endomorphismes commute avec celle du groupe de Galois et en particulier du groupe d'inertie, ce qui entra\^{\i}ne la premi\`ere affirmation. De plus, pour tout endomorphisme $\alpha$, on a 
$$<\alpha v,v'>=<v,\alpha^{\dagger}v'>,$$ o\`u $<\cdot,\cdot>$ d\'esigne l'accouplement de Weil et $\dagger$ l'involution de Rosati associ\'ee \`a la polarisation choisie; ainsi les d\'ecompositions de la repr\'esentation galoisienne d\'ecrites dans  la proposition (\ref{2copies}) induisent celles sur $U_1$ et $U_2$.
\hfill$\Box$

\medskip

L'hypoth\`ese $(T)$ se traduit en disant que $de=\rang U_2=d\sum_{\lambda}\rang U_{i,\lambda}$, ce qui impose $\rang U_{i,\lambda}=1$. Si l'on \'ecrit maintenant la matrice de $\rho(g_I)$ par blocs, on voit que chaque bloc est une transvection. Les arguments de \cite{hall} s'appliquent alors, permettant de montrer que chaque bloc de la repr\'esentation modulo $\ell$ contient dans son  image $\Sp_{2h}(\F_{\ell})$, ce qui ach\`eve la preuve.
\hfill$\Box$

\begin{theo}\label{t2} Soit $s\geq 1$ un entier et soient $A_1,\ldots,A_s$ des vari\'et\'es ab\'eliennes, deux \`a deux non isog\`enes, chaque $A_i$ \'etant pleinement de type Lefschetz, de type I ou II, de dimension relative un entier $h_i$. On a dans ce cas
\[\Hdg\left(\prod_{i=1}^sA_i\right)=\prod_{i=1}^s\Hdg(A_i)=\prod_{i=1}^s\Res_{E_i/\Q}\Sp_{2h_i,E_i}\]
\[\forall \ell,\ H_{\ell}\left(\prod_{i=1}^sA_i\right)=\prod_{i=1}^sH_{\ell}(A_i)=\prod_{i=1}^s\prod_{\lambda_i\,|\,\ell}\Sp_{2h_i,E_{\lambda_i}}\]
o\`u $\lambda_i$ parcourt les places de $E_i$ au dessus de $\ell$.

\end{theo}
\demo Nous expliquons la preuve dans le cadre $\ell$-adique qui s'appuie sur l'article de Lombardo \cite{lombardo} (le cas complexe est plus simple et s'appuie de mani\`ere parall\`ele sur l'article ant\'erieur d'Ichikawa \cite{ichi}). Il s'agit en fait d'appliquer le th\'eor\`eme 4.1 de \cite{lombardo}. Pour cela il nous suffit de v\'erifier que les hypoth\`eses de ce th\'eor\`eme sont satisfaites. Nos hypoth\`eses entra\^inent que pour tout entier $i$, on a $\Hdg(A_i)\times \C_{\ell}=\Sp_{2h_i,\C_{\ell}}^{[E_i:\Q]}$ o\`u $E_i$ est le centre de $\End(A_i)\otimes\Q$. Or les automorphismes de (l'alg\`ebre de Lie de) $\Sp_{2h_i,\C_{\ell}}$ sont int\'erieurs et les automorphismes int\'erieurs pr\'eservent les plus haut poids (cf. remarque 3.8 de \cite{lombardo}), donc le point 3. des hypoth\`eses du th\'eor\`eme 4.1 de \cite{lombardo} est v\'erifi\'e. Les points 1. et 2. sont imm\'ediats dans notre situation, d'o\`u la conclusion. \hfill$\Box$

\medskip

\noindent En particulier ceci nous donne le corollaire suivant :

\begin{cor}\label{c2} Soit $s\geq 1$ un entier et soient $A_1,\ldots,A_s$ des vari\'et\'es ab\'eliennes, deux \`a deux non isog\`enes, chaque $A_i$ \'etant de type I ou II, de dimension relative un entier $h_i$. Pour chaque $i$, on suppose que l'une des  hypoth\`eses suivantes est v\'erifi\'ee : 
\begin{enumerate}
\item l'entier $h_i$ n'est pas dans l'ensemble exceptionnel $\Sigma$ et le centre de $\End(A)\otimes \Q$ est r\'eduit \`a $\Q$,
\item l'entier $h_i$ est \'egal \`a deux ou est impair,
\item la vari\'et\'e ab\'elienne $A_i$ est de type I (resp. II) et poss\`ede une place de mauvaise r\'eduction semi-stable avec dimension torique $e_i$ (resp. avec dimension torique $2e_i$).  
\end{enumerate}
\noindent Sous ces hypoth\`eses on a alors
\[\Hdg\left(\prod_{i=1}^sA_i\right)=\prod_{i=1}^s\Hdg(A_i)=\prod_{i=1}^s\Res_{E_i/\Q}\Sp_{2h_i,E_i}\]
\[\forall \ell,\ H_{\ell}\left(\prod_{i=1}^sA_i\right)=\prod_{i=1}^sH_{\ell}(A_i)=\prod_{i=1}^s\prod_{\lambda_i\,|\,\ell}\Sp_{2h_i,E_{\lambda_i}}\]
o\`u $\lambda_i$ parcourt les places de $E_i$ au dessus de $\ell$.
\end{cor}
\demo Soit $i$ un entier dans l'ensemble $\{1,\ldots,s\}$. Si l'hypoth\`ese $1.$ est v\'erifi\'ee alors par le th\'eor\`eme 5.14 de \cite{pink} (dans le cas de type I) ou par le th\'eor\`eme \ref{t1} ci-dessus (dans le cas de type II) on en d\'eduit que $A$ est pleinement de type Lefschetz, de groupe de Hodge, $\Sp_{2h_i}$. Si l'hypoth\`ese 2. est v\'erifi\'ee, l\`a encore la m\^eme conclusion vaut en appliquant cette fois le th\'eor\`eme A de \cite{BGK} (cf. \cite{lombardo} remarque 2.25 pour le cas de dimension relative 2). Il suffit donc d'appliquer le th\'eor\`eme \ref{t2} ci-dessus pour conclure.\hfill$\Box$

\medskip

\rem On peut bien s\^ur \'enoncer le th\'eor\`eme et le corollaire pr\'ec\'edents en rempla\c{c}ant  $\prod_{i=1}^sA_i$ par $\prod_{i=1}^sA_i^{m_i}$ puisque
$\Hdg\left(\prod_{i=1}^sA_i^{m_i}\right)=\Hdg\left(\prod_{i=1}^sA_i \right)$ et $H_{\ell}\left(\prod_{i=1}^sA_i^{m_i}\right)=H_{\ell}\left(\prod_{i=1}^sA_i\right)$.

\end{document}